\documentclass[10pt]{amsart}
\usepackage{amsmath}
\usepackage[curve,matrix,arrow]{xy}
\title[Recursive Formulas for Characteristic Numbers]
{Recursive Formulas for the Characteristic Numbers\\
of Rational Plane Curves}
\author{Lars Ernstr{\"o}m}
\address{Ohio State University at Mansfield, 1680 Univ. Dr.,
Mansfield, Ohio 44906, USA}
\email{ernstrom@math.ohio-state.edu}
\author{Gary Kennedy}
\address{Ohio State University at Mansfield, 1680 Univ. Dr.,
Mansfield, Ohio 44906, USA}
\email{kennedy@math.ohio-state.edu}

\newcommand{\abstracttext}{We derive recursive equations for the
characteristic numbers of rational nodal plane curves with at most one cusp,
subject to point conditions, tangent conditions and flag conditions,
developing techniques akin to quantum cohomology on a moduli space of
stable lifts.}

\numberwithin{equation}{section}
\newtheorem{thm}[equation]{Theorem}
\newtheorem{prop}[equation]{Proposition}
\newtheorem{cor}[equation]{Corollary}

\theoremstyle{definition}

\theoremstyle{remark}

\newcommand{\pp}{{\mathbf P^2}} 
\newcommand{\pl}{{\mathbf P^1}} 
\newcommand{\bp}{\mathbf P} 
\newcommand{\dpp}{\mathbf {\check P}^{2}} 
\newcommand{\q}{\mathbf Q} 
\newcommand{\hd}{\check{h}}
\newcommand{\OO}{\mathcal O}
\newcommand{\lmoduli}{{\overline M}^1_{0,n}(\pp,d)}
\newcommand{\lamoduli}[2]{{\overline M}^1_{0,#1}(\pp,#2)}
\newcommand{\imoduli}{{\overline M}_{0,n}(I,(d,2d-2))}
\newcommand{\icmoduli}{{\overline M}_{0,n}(I,(0,0))}
\newcommand{\ppmoduli}{{\overline M}_{0,n}(\pp,d)}

\newcommand{\stablecurve}{{\overline M}_{0,n}}
\newcommand{\cusplocus}{{\overline C}_{0,n,\{\star\}}(\pp,d)}
\newcommand{\lcusplocus}{{\overline C}^1_{0,n,\{\star\}}(\pp,d)}
\newcommand{\ldivisor}{D\left(A_1,A_2,(d_1,c_1),(d_2,c_2)\right)}

\newcommand{\codim}{\operatorname{codim}}
\newcommand{\xmoduli}{{\overline M}_{0,n}(X,\varphi)}
\newcommand{\ymoduli}{{\overline M}_{0,n}(Y,f_*(\varphi))}
\newcommand{\stablefour}{{\overline M}_{0,4}}

\newcommand{\amoduli}[2]{{\overline M}_{0,#1}(X,#2)}

\newcommand{\calP}{{\mathcal P}}
\newcommand{\calN}{{\mathcal N}}
\newcommand{\calC}{{\mathcal C}}
\newcommand{\calF}{{\mathcal F}}
\newcommand{\calR}{{\mathcal R}}
\newcommand{\calL}{{\mathcal L}}
\newcommand{\calE}{{\mathcal E}}
\newcommand{\calK}{{\mathcal K}}
\newcommand{\calG}{{\mathcal G}}
\begin{document}
\begin{abstract} \abstracttext \end{abstract} \maketitle
%
%
%
\section{Introduction}
\par
Although of recent vintage, Kontsevich's recursive formula $$
N_d=\sum_{d_1+d_2=d}N_{d_1}N_{d_2}\left[d_1^2d_2^2\binom{3d-
4}{3d_1-2}- d_1^3d_2\binom{3d-4}{3d_1-1}\right]
$$
for the number $N_d$ of rational plane curves of degree $d$ through $3d-1$
general points is already a celebrated result. The formula is wonderfully
simple: it determines all such characteristic numbers, beginning with the
triviality $N_1 = 1$.
\par
In his proof, Kontsevich uses a compactification $\ppmoduli$ of the space of
$n$-pointed maps from $\pl$ to the projective plane, called the space of {\em
stable maps}. The key fact about these spaces is that (for $n \geq 4$)
they carry three linearly equivalent divisors, called the {\em special boundary
divisors}, which are the inverse images of the special points $0$, $1$, and
$\infty$, under a natural map from $\ppmoduli$ to the space of stable
rational curves with four distinguished points. (As is well-known, the latter
space is isomorphic to $\pl$.) This linear equivalence readily yields
Kontsevich's identity. Furthermore, each divisor is a union of components,
each of which is isomorphic to a fiber product
$$
{\overline M}_{0,n_1+1}(\pp,d_1)
\times_{\pp}
{\overline M}_{0,n_2+1}(\pp,d_2);
$$
as a consequence, one can interpret Kontsevich's formula as the assertion that
a certain new product on the cohomology group of $\pp$ (tensored with a
certain power series ring) is associative. The resulting associative ring is
called
the {\em quantum cohomology}.
\par
In this paper we show that Kontsevich's compactifications can also be used to
derive recursive formulas for these other classical characteristic numbers:
\begin{align*}
N_d(a,b,c)=&\text{ the number of rational plane curves of degree $d$
through $a$ general points, tangent to}\\ &\text{$b$ general lines, and
tangent to $c$ general lines at a specified general point on each line}\\
&\quad\text{(where $a+b+2c=3d-1$).}
\end{align*}
\par
In examining questions of tangency, it is natural to work with the incidence
correspondence $I$ of points and lines in $\pp$. However, the space of stable
maps to $I$ is too large; instead we need to consider only those stable maps
which can be lifted from maps to $\pp$, and degenerations of such maps.
Thus we will define a space $\lmoduli$ of {\em stable lifts}. The superscript
$1$ indicates that $I$ parametrizes first-order jets of curves in $\pp$;
similar definitions for higher order jet bundles will be considered in another
paper. As a subspace of
the space of stable maps to the incidence correspondence, $\lmoduli$ inherits
the special boundary divisors. A chief part of our project is to describe the
components of these divisors. We show, for example, that for many
components the general point represents a map from a curve with {\em
three} components, with the central component mapping to a fiber of $I$ over
$\pp$.
\par From the linear equivalence of the special boundary divisors, we
extract many recursive formulas. We show, in an admittedly {\it ad hoc}
fashion, that these formulas suffice to determine all the characteristic
numbers $N_d(a,b,c)$, starting with the four trivial cases in degree $1$.
 At the same time we obtain recursive formulas which determine, in all cases,
three other types of characteristic numbers:
\begin{align*}
C_d(a,b,c;1)=&\text{ the number of rational plane curves of degree $d$, with
one cusp, through $a$ points,}\\
&\text{tangent to $b$ lines, and tangent to $c$ lines at $c$ specified
points}\\
&\quad\text{(where $a+b+2c=3d-2$);}\\ C_d(a,b,c;h)=&\text{ the number of
such curves having the cusp on a specified line}\\
&\quad\text{(where $a+b+2c=3d-3$);}\\
C_d(a,b,c;h^2)=&\text{ the number of such curves having the cusp at a
specified point}\\
&\quad\text{(where $a+b+2c=3d-4$).}
\end{align*}
In fact, it is easier to determine these numbers than to avoid them, since the
special boundary divisors on the space of stable lifts include, in addition to
the
components already described, other components corresponding to cuspidal
rational curves.
\par
Because of the occurrence of stable maps with automorphism, it is most
convenient to utilize the framework of algebraic stacks and the corresponding
intersection theory (with rational coefficients) as defined by Vistoli
\cite{Vistoli}. For automorphism-free stable maps the intersection theory of
the stack coincides with the usual intersection theory of the corresponding
moduli scheme. However, on the stack $\lmoduli$ of
stable lifts there are components of the special boundary divisors for which
the general member has automorphism group of order two. Therefore we will
encounter the fractional intersection number $1/2$. Another possibility would
be to work with equivariant Chow groups as defined by Edidin and Graham
\cite{EdidinGraham}.
\par
Here is an outline of the paper.
In \S\ref{srmaps} we recall Kontsevich's notion of stable maps, the
definition of the Gromov-Witten invariants, and the basic linear equivalence
which leads to Kontsevich's recursive formula. As we have said, this linear
equivalence implies the associativity of the quantum product,
but---since in our project we do not use any analogous ``higher-order''
product---we refrain from discussing this point.
In \S\ref{sectionlifts} we introduce spaces of stable lifts for $\pp$, by
associating to a general stable map its lift to the projectivized tangent
bundle of  $\pp$.
In \S\ref{sectionHOGW} we use these spaces to define first-order Gromov-Witten
invariants, and discuss their significance in enumerative geometry.
In \S\ref{sectiondivisors} we describe the structure of the special boundary
divisors on the space of stable lifts.
In \S\ref{sectionpotentials} we use the
physicists' method of quantum potentials to write down generating functions
for the Gromov-Witten invariants, and derive the basic identity
(\ref{ijklrel}), analogous to the associativity identity for the quantum
product.
In \S\ref{sectionrecursion} we derive, from selected cases of the basic
identity, recursive equations for the characteristic numbers;
we also present tables of characteristic numbers through degree $5$.
\par
We believe that these methods will work for a large class of surfaces. We also
believe that they will enable us to calculate the higher-order characteristic
numbers of curves, defined by conditions of higher-order tangency to lines or
other specified curves. Here the appropriate parameter spaces are, we believe,
those originally defined by Semple (later investigated in \cite{Collino} and
\cite{CKtrip}). We intend to pursue this matter with Susan Colley. A theory
for more general varieties and higher genus seems possible as a consequence
of the results of Behrend \cite{Behrend} and those of Li and Tian
\cite{LiTian2}.
However, the enumerative significance of the Gromov-Witten invariants is,
in this more general set-up, not obvious.
\par
In recent papers, Pandharipande has constructed an algorithm for computing the
numbers $N_d(a,b,0)$ \cite{Pandharipande2}, and Francesco and Itzykson have
found a recursive relation that determines the same numbers, with the
restriction $a\geq 3$, given characteristic numbers of lower degree
\cite{FrancescoItzykson}.
\par
This project was inspired by W.~Fulton's lecture series at the AMS summer
institute at Santa Cruz 1995.
We have benefited greatly from discussions with R.~Pandharipande,
and from his preprints. We also wish to thank P.~Aluffi, C.~Ban,
S.~Colley, B.~Crauder,
D.~Edidin, S.~Kleiman, D.~Laksov, L.~McEwan, and S.~Yokura for their help
during the preparation of this paper.
\par
%
%
%
\section{Stable maps and Gromov-Witten invariants}\label{srmaps}
As general references for the material in this section we suggest
\cite{BehrendManin}, \cite{Fulton2}, \cite{Kontsevich}, \cite{KontsevichManin},
\cite{LiTian} and  \cite{Pandharipande}. In this section we will utilize the
intersection theory of the coarse moduli space of stable maps. Generic members
of the special boundary divisors of the space of stable
maps have no automorphisms. Therefore, for our purposes, the intersection
theory of the stack of stable maps is the same as that of the coarse moduli
space. Throughout this paper, we will be working over an algebraically closed
field of characteristic zero. Given a variety or stack $X$, homology $A_*(X)$
and cohomology $A^*(X)$ will denote the rational equivalence groups with
coefficients in $\q$.
\par
Suppose that $C$ is a connected reduced curve of arithmetic genus zero
whose singularities are at worst nodes. Then $C$ must be a tree of $\pl$'s.
We will call the nodes {\em points of attachment}. Suppose that on $C$ we
have $n$ distinct points $p_1,\dots,p_n$, none of which is a point of
attachment; we call these points {\em markings}. A {\em special point} is a
point of attachment or a marking.
\par
Now suppose that $X$ is a nonsingular projective variety. Following
Kontsevich, we define an {\em $n$-pointed stable (genus $0$) map to $X$} to be
a
map $\mu\colon C \to X$ from such a curve, subject to the following
condition: if the restriction of $\mu$ to a component is constant, then that
component contains at least three special points. (Since we never consider
stable maps from curves of higher arithmetic genus, we will generally omit the
parenthesized phrase.) A {\em family of stable maps} consists of a flat, proper
map $\pi\colon \mathcal C\to S$, a
map $\mu\colon\mathcal C\to X$, and $n$ sections
$\{p_t\}_{t=1,\dots,n}$ of $\pi$, such that, for each geometric fiber
$\mathcal C(s)$
of $\pi$, the restriction of $\mu$ to this fiber, with the markings $p_t(s)$,
is
an $n$-pointed stable map. A {\em morphism} over $S$ from one
such family $(\pi\colon \mathcal C \to S,\mu,\{p_t\})$ to another
$(\pi^\prime\colon \mathcal C^\prime \to
S,\mu^\prime,\{p^\prime_t\})$
is a morphism $\tau\colon \mathcal C\to \mathcal C^\prime$ over $S$
such that $\mu=\mu^\prime\circ \tau$ and $p_t\circ \tau=p_t^\prime$.
\par
Suppose that $\varphi$ is a specified class in $A_1(X)$. Kontsevich
 showed that there is a coarse moduli space $\xmoduli$ for
isomorphism classes of stable maps whose image in $X$ represents
$\varphi$ \cite[2.4.1]{KontsevichManin}, \cite[\S1, p. 336]{Kontsevich},
\cite[Theorems 1 and 2]{Pandharipande}.
We call $\xmoduli$ the {\em space of stable (genus $0$) maps}. There are $n$
{\em evaluation maps} $e_t:\xmoduli\to X$ which associate to the point
representing a stable map $\mu$ the images $\mu(p_t)$ of the
markings. If $f:X \to Y$ is a morphism, then there is a morphism
$Mf:\xmoduli \to \ymoduli$ which associates to the point representing
$\mu$ the point representing the stable map obtained from $f \circ \mu$ by
contracting any components which have become unstable; if $f_*(\varphi)=0$, we
must assume that $n \geq 3$. In particular, under this hypothesis
there is a ``forgetful'' morphism from $\xmoduli$ to $\stablecurve$, the
space of stable rational $n$-pointed curves. \par If $n \geq 4$ then we can
compose the forgetful morphism with another sort of forgetful morphism
$\stablecurve \to \stablefour$, this time forgetting about all markings except
the first four. The space $\stablefour$ is isomorphic to $\pl$. It has a
distinguished point $P(12\mid 34)$ representing the two-component curve
having the first two markings on one component and the latter two on the
other; similarly there are two other distinguished points $P(13\mid 24)$ and
$P(14\mid 23)$. Their inverse images on $\xmoduli$ are three linearly
equivalent divisors $D(12\mid 34)$, $D(13\mid 24)$, and $D(14\mid 23)$,
called the {\em special boundary divisors}.
\par We allow the index set to be an arbitrary finite
set $A$ rather than just $\{1,\dots,n\}$; we then write
$\amoduli{A}{\varphi}$ for the moduli space. Suppose that $A_1 \cup A_2$
is a partition of $\{1,\dots,n\}$ and that $\varphi_1$, $\varphi_2$ are
classes in $A_1(X)$ whose sum is $\varphi$. Suppose that $\{\star\}$ is a
single-element set. Then the fiber product $$
D(A_1,A_2,\varphi_1,\varphi_2) =
\amoduli{A_1\cup\{\star\}}{\varphi_1} \times_X
\amoduli{A_2\cup\{\star\}}{\varphi_2} $$ is naturally a subspace of
$\amoduli{A}{\varphi}$; the typical point represents a map from a curve
with two components, with the point of attachment corresponding to the
point indexed by $\star$. In Figure 2.1 the attachment point and the image
thereof are marked $\bullet$.
$$
\xy
0;<0.5cm,0cm>:
(-1,-1);(1,1)**\dir{-};
(-1,1);(1,-1)**\dir{-};
(0,0)*{\bullet};
(1.5,0);(3.5,0)**\dir{-};
?>*\dir{>};
(4,-0.5);(6.5,0.75)**\dir{-};
(4,1);(6,1)**\crv{(5,-3)};
(4.41,-0.3)*{\bullet}
\endxy
$$
\begin{center}
{\bf Figure 2.1} Stable map.
\end{center}
\par
\bigskip
If $X$ is a product of homogeneous
spaces (for example, if $X$ is a projective space or a flag variety) then the
dimension of $\xmoduli$ (if it is nonempty) is
$\dim X + \int_\varphi c_1(T_X) + n - 3$ \cite{Kontsevich}, \cite[Proposition
1]{Fulton2}. Furthermore the divisor
$D(12\mid 34)$ is the sum
\begin{equation}\label{donetwo} D(12\mid
34)=\sum D(A_1,A_2,\varphi_1,\varphi_2) \end{equation}
over all
partitions in which $1$ and $2$ belong to $A_1$, and $3$ and $4$ belong to
$A_2$.
\par
Suppose that
$\gamma_1,\dots,\gamma_n$ are elements of $A^*X \otimes \q$. Then
the {\em Gromov-Witten invariant} associated to these classes and to
$\varphi$ is the number
$$
N_{\varphi}(\gamma_1\cdots\gamma_n)
=\int e_1^*(\gamma_1) \cup \dots \cup e_n^*(\gamma_n)
\cap [\xmoduli],
$$
the degree of the top-dimensional component. If the classes
$\gamma_1,\dots,\gamma_n$ are homogeneous, then the Gromov-Witten
invariant is nonzero only if the sum of their codimensions is the dimension
of the moduli space.
(If the moduli space is empty, we declare the Gromov-Witten invariant to be
zero.)
\par
We now restrict our attention to the case $X=\pp$.
The class $\varphi$ of a curve must be some multiple $dh$ of the class $h$ of
a line. We will write $\ppmoduli$ for the space of $n$-pointed stable maps; its
dimension is $3d-1+n$ (unless $d=0$ and $n \leq 2$, in which case
the space is empty). We define
$$
N_d=N_{dh}(h^2 \cdots h^2).
$$
By standard transversality arguments, one can show that $N_d$ is the
number of rational plane curves of degree $d$ through $3d-1$ general points.
\par
The linear equivalence of the special boundary divisors
$D(12\mid 34)$ and $D(13\mid 24)$ implies, for each
choice of $\varphi$ and of $\gamma_1,\dots,\gamma_n$,
the numerical equality
\begin{equation}
\label{numeq}
\int e_1^*(\gamma_1) \cup \dots \cup e_n^*(\gamma_n)
\cap [D(12\mid 34)] =
\int e_1^*(\gamma_1) \cup \dots \cup e_n^*(\gamma_n)
\cap [D(13\mid 24)].
\end{equation}
By (\ref{donetwo}) and the similar decomposition of the divisor $D(13\mid
24)$, this is an equation among various values of $N_d$. For example, if
$n=3d+1$, $\varphi=dh$, $\gamma_1=\gamma_2=h$, and
$\gamma_t=h^2$ for $t \geq 3$, then (\ref{numeq}) is Kontsevich's
recursive formula!
\par
Equation (\ref{numeq}) can also be interpreted as an equation among formal
power series. Let $T_0=1$, $T_1=h$, and $T_2=h^2$, so that an arbitrary class
$\gamma$ is a linear combination $y_0 T_0 + y_1 T_1 + y_2 T_2$. Define
the {\em potential} by
\begin{equation}
\calP=\sum_{n \geq 0}\frac{1}{n!}
\sum_{\varphi}
N_{\varphi}(\gamma^n),
\end{equation}
which one can show is equal to
$$
\frac{1}{2}(y_0^2y_2+y_0y_1^2)
+ \sum_{d \geq 1} N_d e^{dy_1}\frac{y_2^{3d-1}}{(3d-1)!}.
$$
The first term, called the {\em classical potential}, encodes the intersection
product on $\pp$; the second term is called the {\em quantum potential}.
Using these potentials, we can translate the special case
\begin{align*}
\int e_1^*(T_i) & \cup e_2^*(T_j) \cup e_3^*(T_k) \cup e_4^*(T_l)
\cup e_5^*(\gamma) \cup \dots \cup e_n^*(\gamma)
\cap [D(12\mid 34)] \\
&=
\int e_1^*(T_i) \cup e_2^*(T_j) \cup e_3^*(T_k) \cup e_4^*(T_l)
\cup e_5^*(\gamma) \cup \dots \cup e_n^*(\gamma)
\cap [D(13\mid 24)]
\end{align*}
of (\ref{numeq}) into the partial differential equation
\begin{equation}
\sum_{s=0}^2 \frac{\partial^3\calP}{\partial y_i\partial y_j\partial y_s}
\frac{\partial^3\calP}{\partial y_k\partial y_l\partial y_{2-s}} =
\sum_{s=0}^2 \frac{\partial^3\calP}{\partial y_i\partial y_k\partial y_s}
\frac{\partial^3\calP}{\partial y_j\partial y_l\partial y_{2-s}}.
\end{equation}
The case $i=j=2$, $k=l=1$ is equivalent to Kontsevich's formula.
%
%
%
%
\section{Lifts of stable maps}\label{sectionlifts}
Denote by $I=\bp(T_{\pp})$ the projectivized tangent
bundle of the projective plane; it is the variety of point-line incidence in
$\pp\times\dpp$.
Suppose that $\mu \colon (\pl,p_1,\dots,p_n) \to \pp$ is a nonconstant
map with $n$ distinct marked points.
For a general point $x$ in $\pl$
define $\tilde\mu\colon\pl\to I$ by
$\tilde\mu(x)=(\mu(x),\mu^\prime(x))$, where $\mu^\prime(x)$ is
the tangent line to $\mu(\pl)$ at $\mu(x)$. The construction extends to all
of $\pl$: even at points where the map $\mu$ is singular there is a unique
tangent direction associated to the branch of $\mu(\pl)$ at $x$.
The map $\tilde\mu\colon (\pl,p_1,\dots,p_n)\to I$ is called the {\em
strict lift} of $\mu$. (By symmetry, we may also define the strict lift of a
nonconstant map to $\dpp$.)
$$
\xy
0;<0.5cm,0cm>:
(-10,0)*{\pl};
(-8.5,0);(-6.5,0)**\dir{-};
?>*\dir{>};
(-7.5,-0.5)*{\mu};
(-5,0)*{\pp};
(-8.5,1.5);(-6.5,3.5)**\dir{-};
?>*\dir{>};
(-5,5)*{I};
(-5,3.5);(-5,1.5)**\dir{-};
?>*\dir{>};
(-8,3)*{\tilde\mu};
(-1,1);(1,-1)**\dir{-};
(1.5,0);(3.5,0)**\dir{-};
?>*\dir{>};
(4,0.5);(5.5,0)**\crv{(8,-0.5)&(4,-1.5)};
(5.5,0);(6,0.5)**\dir{-};
(1.5,1.5);(3.5,3.5)**\dir{-};
?>*\dir{>};
(5,3.5);(5,1.5)**\dir{-};
?>*\dir{>};
(4,5.5);(5.4,4.9)**\crv{(8,4.5)&(4,3.5)};
(5.55,5.05);(6,5.5)**\dir{-};
\endxy
$$
\begin{center}
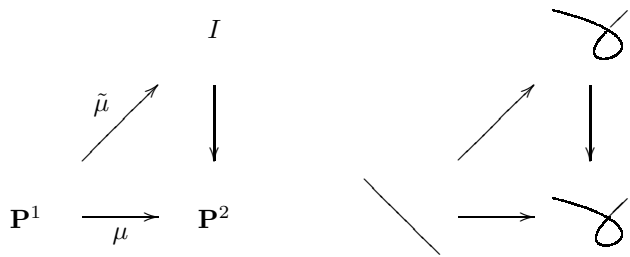
{\bf Figure 3.1} The strict lift of a map from $\pl$ to $\pp$.
\end{center}
\par
\bigskip
Strict lifts do not behave nicely in families of maps. For example if $\mu$ is
an immersion of degree $d$ and its image is a nodal curve, then the number
of nodes is $\delta=(d-1)(d-2)/2$ and the class of the curve is
$$\check{d}=d(d-1)-2\delta=2d-2.$$
Thus the homology class of $\tilde \mu(\pl)$
in $I$ is $d$ times the strict lift of a line plus $2d-2$ times the strict lift
of a
dual line.
We will say that it has {\em bidegree} $(d,2d-2)$. If instead $\mu(\pl)$ is a
rational curve with $(d-1)(d-2)/2-1$ nodes and one cusp then the bidegree of
$\tilde\mu(\pl)$ is
$(d,2d-3)$. Thus in a family degenerating a node to a cusp, the strict lifts of
the
members in the family do not piece together.
\par
Suppose that we have a family of stable maps $(\pi\colon \mathcal C \to
S,\mu\colon\mathcal C\to I,\{p_t\})$ whose general member is the strict
lift of an immersion $\pl\to \pp$. Then we say that the members of this
family are {\em stable lifts} of the corresponding members of the family
of maps to $\pp$. Note that a stable lift of a map of degree $d$ has bidegree
$(d,2d-2)$; for example (supposing there are no markings) the stable lift of a
$d$-fold branched cover of a line consists of the strict lift together with
maps to the fibers of $I$ over the $2d-2$ branch points. In Figure 3.2 we
illustrate the case $d=3$. Ramification and branch points are marked by
$\times$.
\par
$$
\xy
0;<0.5cm,0cm>:
(-1,0);(1,0)**\dir{-};
(-0.6,-1);(-0.6,1)**\dir{-};
(-0.6,0)*{\bullet};
(-0.6,0)*{\times};
(-0.2,-1);(-0.2,1)**\dir{-};
(-0.2,0)*{\bullet};
(-0.2,0)*{\times};
(0.2,-1);(0.2,1)**\dir{-};
(0.2,0)*{\bullet};
(0.2,0)*{\times};
(0.6,-1);(0.6,1)**\dir{-};
(0.6,0)*{\bullet};
(0.6,0)*{\times};
(1.5,0);(3.5,0)**\dir{-};
?>*\dir{>};
(4,0);(6,0)**\dir{-};
(4.4,-1);(4.4,1)**\dir{-};
(4.4,0)*{\bullet};
(4.8,-1);(4.8,1)**\dir{-};
(4.8,0)*{\bullet};
(5.2,-1);(5.2,1)**\dir{-};
(5.2,0)*{\bullet};
(5.6,-1);(5.6,1)**\dir{-};
(5.6,0)*{\bullet};
(4.4,0)*{\times};
(4.8,0)*{\times};
(5.2,0)*{\times};
(5.6,0)*{\times};
(0,-1.5);(0,-3.5)**\dir{-};
?>*\dir{>};
(5,-1.5);(5,-3.5)**\dir{-};
?>*\dir{>};
(-1,-5);(1,-5)**\dir{-};
(-0.6,-5)*{\times};
(-0.2,-5)*{\times};
(0.2,-5)*{\times};
(0.6,-5)*{\times};
(1.5,-5);(3.5,-5)**\dir{-};
?>*\dir{>};
(4,-5);(6,-5)**\dir{-};
(4.4,-5)*{\times};
(4.8,-5)*{\times};
(5.2,-5)*{\times};
(5.6,-5)*{\times};
\endxy
$$
\begin{center}
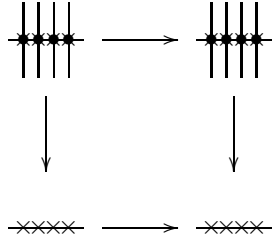
{\bf Figure 3.2} The stable lift of a three-cover of a line in $\pp$ with
ramification at four points.
\end{center}
\par
\bigskip
For simplicity of notation, we will write $\imoduli$ for the stack of stable
(genus $0$) maps representing the class of a curve of bidegree $(d,2d-2)$. Let
$\lmoduli$ be the closed substack of $\imoduli$ which represents stable lifts;
we
call it the {\em stack of stable lifts}. We continue to use the notations
$e_1,\dots,e_n$ for the evaluation maps $\lmoduli \to I$.
Note that the map $\lmoduli \to \ppmoduli$ (inclusion followed by
projection) is a birational morphism, whose inverse is the {\em lifting map}
$\lambda$ which associates to each immersion its strict lift.
\begin{prop}\label{stableliftd1d2}
Let $C$ be a curve consisting of two irreducible component, each isomorphic to
$\pl$ and attached at a single point $p$.
Let $\mu\colon C \to \pp$ be a map defined by mapping one component of $C$
to an irreducible curve of degree $d_1\geq 1$ and mapping the other
component to an irreducible curve of degree $d_2\geq 1$, so that $\mu$ is
represented by a point of ${\overline M}_{0,0}(\pp,d_1+d_2)$.
Assume that both maps are immersions and that they are transverse at $p$.
\par
Then the stable lift of $\mu$ is unique. It is a map $\tilde\mu$ from a curve
$\tilde C$ consisting of three components, each isomorphic to $\pl$. On the
cental component $\tilde\mu$ is a map of degree $2$ to the pencil of tangent
directions at $p$, with the points of attachment mapping to the tangent
directions of $C$; the map is ramified at these attachment points. On each
peripheral component $\tilde\mu$ is the strict lift of a component of $\mu$.
\end{prop}
\par
See Figure 3.3.
\par
$$
\xy
0;<0.5cm,0cm>:
(-1,-1);(1,1)**\dir{-};
(-1,1);(1,-1)**\dir{-};
(0,0)*{\bullet};
(1.5,0);(3.5,0)**\dir{-};
?>*\dir{>};
(4,-0.5);(6.5,0.75)**\dir{-};
(4,1);(6,1)**\crv{(5,-3)};
(4.41,-0.3)*{\bullet};
(0,6);(0,4)**\dir{-};
(0,5.333)*{\bullet};
(0,5.333)*{\times};
(0,4.667)*{\bullet};
(0,4.667)*{\times};
(-1,5.333);(1,5.333)**\dir{-};
(-1,4.667);(1,4.667)**\dir{-};
(1.5,5);(3.5,5)**\dir{-};
?>*\dir{>};
(4.41,6);(4.41,4)**\dir{-};
(4,4.667);(6,4.667)**\dir{-};
(4.41,4.667)*{\bullet};
(4.41,4.667)*{\times};
(4,6);(6,6)**\crv{(6,5)};
(4.41,5.78)*{\bullet};
(4.41,5.78)*{\times};
(0,3.5);(0,1.5)**\dir{-};
?>*\dir{>};
(5,3.5);(5,1.5)**\dir{-};
?>*\dir{>};
\endxy
$$
\begin{center}
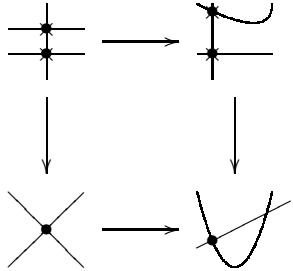
{\bf Figure 3.3} Stable lift of Proposition \ref{stableliftd1d2}.
\end{center}
\par
\bigskip
\begin{proof}
The versal deformation theory of a plane node is given by the local equation
$xy=\epsilon$. Thus any family degenerating to $\mu$ will be the pullback of a
family which, in appropriate local coordinates, has total space isomorphic to a
neighborhood of the origin, with the fibers of the family being the curves
$xy=\epsilon$ and with the map being the identity. Locally around the origin,
$I$ is a trivial $\pl$-bundle; in coordinates the family of stable lifts is
given by
$$ (x,y)\mapsto \bigg((x,y),[x:-y]\bigg).$$
\par
According to \cite[section 3.2]{Pandharipande}, there is a unique extension of
the family of lifts, obtained by blowing up the points of indeterminacy,
possibly after a base change. We choose to make the base change which has the
effect of replacing our original family by the family $xy=\epsilon^2$. (In
fact, one can show that this base change is unavoidable.) Then we blow up the
origin on this singular surface. In the resulting family, the fiber over
$\epsilon=0$ is reduced and has a new component mapping to the pencil of
directions at the origin; the map, which is of degree 2, is ramified at the
points of attachment to the strict transforms of the two original components.
\end{proof}
But in general a stable map may have more than one stable lift.
Figure 3.4 illustrates a stable lift of a degree three cover of a line with
double ramification at two points. The two maps into the fibers over the branch
points are of degree two, ramified at the attachment point and at one other
point. The latter point may appear anywhere on the fiber. Thus the stable lift
in this situation is not unique; there is a two dimensional family of
lifts corresponding to a single stable map.
$$
\xy
0;<0.5cm,0cm>:
(-1,0.333);(1,0.333)**\dir{-};
(-0.333,-1);(-0.333,1)**\dir{-};
(-0.333,0.333)*{\bullet};
(-0.333,0.333)*{\times};
(-0.333,-0.333)*{\times};
(0.333,-1);(0.333,1)**\dir{-};
(0.333,0.333)*{\bullet};
(0.333,0.333)*{\times};
(0.333,-0.333)*{\times};
(1.5,0);(3.5,0)**\dir{-};
?>*\dir{>};
(4,0.333);(6,0.333)**\dir{-};
(4.667,-1);(4.667,1)**\dir{-};
(4.667,0.333)*{\bullet};
(4.667,0.333)*{\times};
(4.667,-0.333)*{\times};
(5.333,-1);(5.333,1)**\dir{-};
(5.333,0.333)*{\bullet};
(5.333,0.333)*{\times};
(5.333,-0.333)*{\times};
(0,-1.5);(0,-3.5)**\dir{-};
?>*\dir{>};
(5,-1.5);(5,-3.5)**\dir{-};
?>*\dir{>};
(-1,-5);(1,-5)**\dir{-};
(-0.333,-5)*{\times};
(0.333,-5)*{\times};
(1.5,-5);(3.5,-5)**\dir{-};
?>*\dir{>};
(4,-5);(6,-5)**\dir{-};
(4.667,-5)*{\times};
(5.333,-5)*{\times};
\endxy
$$
\begin{center}
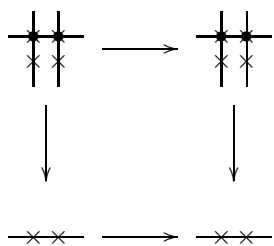
{\bf Figure 3.4} A stable lift of a map with double ramification is not
unique.
\end{center}
\par
\bigskip
\begin{prop}\label{liftclassif} Let $\tilde\mu\colon \tilde C\to I$ be an
$n$-pointed stable lift. Let $\mu\colon C \to \pp$ be the stable map
obtained by composing $\mu$ with the projection $I\to \pp$, forgetting
about all markings, and (if necessary) contracting components which have
become unstable. Then, for each irreducible component $\tilde C_i$ of $\tilde
C$, the
restricted map $\tilde \mu\vert_{\tilde C_i}\colon \tilde C_i\to I$ is one of
the following:
\begin{trivlist}
\item[(1)] the strict lift of a map obtained by restricting $\mu$ to a
component of $C$,
\item[(2)] a map into a fiber of $I$ over some point $x$ in $\pp$, where $x$ is
either the image of an attachment point of $C$ or a singularity of the
restriction of $\mu$ to some component of $C$,
\item[(3)] a constant map.
\end{trivlist}
\end{prop}
\par
\begin{proof}
We may assume that $n=0$. It then suffices to prove that each component is of
type (1) or (2).
\par
By definition $\tilde\mu$ is a member of a family of stable maps from
$\tilde\calC$ to $I$ over a base $S$, where the generic member of the family is
the strict lift of a map to $\pp$. We may and will assume that the base $S$ is
a nonsingular curve. Consider the family of maps $\calC$ to $\pp$ over the same
base $S$ obtained by composing with the projection $I\to \pp$, forgetting about
all markings, and (if necessary) contracting components which have become
unstable.
We will next reconstruct the family of maps to $I$, and conclude
the facts about the special member $\tilde\mu$.
Consider the open set $U$ of $\calC$ over $S$ that is the complement of the
singular points (the nodes or attachment points) of curves in the family.
Then there is a map
$$T_{U/S}\to T_\pp$$ from the relative tangent bundle of $U$
over $S$ to the tangent bundle of $\pp$.
Composing with the map
$$T_\pp - \text{zero section}\to \bp(T_\pp)=I$$
we get a rational map from $\bp(T_{U/S})=U$ and hence from $\calC$ to $I$,
which is indeterminate exactly at the attachment points and singularities of
the map $\mu$ on components of the special member $C$.
For generic members of $\calC$ the rational map is clearly the lifting map.
Finally, by \cite[Prop.3.3]{Pandharipande} the rational map may be extended by
blowing up the points of indeterminacy. Furthermore by
\cite[Prop.3.2]{Pandharipande} the extension is unique and thus we recover the
map $\tilde\mu$ as the special member after sufficient blowing-up.
The result now follows from what we have concluded about the points of
indeterminacy.
\end{proof}
%
%
%
%
\par
\section{Characteristic numbers}\label{sectionHOGW}
\par
The rational cohomology of $I$ is given by
$$
A^*(I)=\q[h,\hd]/(h^3,\hd^3,h^2+\hd^2-h\hd), $$
where $h$ is the first Chern class of the pullback of the line bundle
$\OO_{\pp}(1)$ on $\pp$
and $\hd$ is the pullback of the line bundle $\OO_{\dpp}(1)$ on $\dpp$.
Note that the fundamental class of the strict lift of a line is $\hd^2$,
and that the class of the strict lift of a dual line is $h^2$.
We fix the following basis for $A^*(I)$:
$$
\{T_0,T_1,T_2,T_3,T_4,T_5\}=\{1,h,h^2,\hd,\hd^2,h^2\hd\}.
$$
With respect to this basis the fundamental class of the diagonal $\Delta$
in $I \times I$ has the simple decomposition
\begin{equation}
\label{diag}
[\Delta] = \sum_{s=0}^5 [T_s] \times [T_{5-s}].
\end{equation}
\par
Suppose that $d$ is a positive integer, and that
$\gamma_1,\dots,\gamma_n$ are elements of $A^*(I) \otimes \q$.
We define the {\em first-order Gromov-Witten invariant} by
$$
N_d(\gamma_1\cdots\gamma_n)=\int e_1^*(\gamma_1) \cup \dots \cup
e_n^*(\gamma_n)
\cap [\lmoduli].
$$
Suppose that $a$ of the $\gamma_t$'s equal the class $h^2$, that $b$ of
them equal the class $\hd^2$, and that the remaining $c$ of them equal
$h^2\hd$, where $a+b+2c=3d-1$. In this case we denote the Gromov-Witten
invariant by $N_d(a,b,c)$.
\par
\begin{thm}\label{gwenumsign2} $N_d(a,b,c)$ is the number of rational plane
curves of degree $d$ through $a$ general points, tangent to $b$ general
lines, and tangent to $c$ general lines at a specified general point on
each line.
\end{thm}
\par
(We will say that such a curve is {\em incident to $c$ specified flags.})
\begin{proof}
Associated to each class $e_t^*h$ and $e_t^*\hd$ there is a complete
basepoint-free linear system parametrized by lines in $\pp$ or by lines in
$\dpp$. Applying the Kleiman-Bertini Theorem repeatedly
\cite[Corollary 5, p.291]{Kleiman4}, we find that the intersection of general
members of the linear systems is regular away from the singular locus of
$\lmoduli$. By a dimension count it is a set of reduced points of
$\lmoduli$. By the same dimension count, each point of the intersection
corresponds to the strict lift of an immersion $\pl \to \pp$. Thus no point
of the singular locus is in the intersection.
\end{proof}
The following Proposition is analogous to \cite[2.2.3 and 2.2.4]
{KontsevichManin} and \cite[p. 9]{Fulton2}.
\begin{prop}\label{GWprop}
\begin{trivlist}
\par
\item[(1)] For all $d\geq 1$ and all $\gamma_1,\dots,\gamma_{n-1}$,
$$
N_d(\gamma_1\cdots\gamma_n\cdot 1)=0. $$
\item[(2)] If $\gamma_n$ is a divisor class in
$A^1(I)$, then for all $d\geq 1$ and all $\gamma_1,\dots,\gamma_{n-1}$,
$$
N_d(\gamma_1\cdots\gamma_n)=
N_d(\gamma_1\cdots\gamma_{n-1})\int\gamma_n\cap [C], $$
where $[C]=d \hd^2 +(2d-2) h^2$ is the class of the strict lift of a
rational nodal plane curve of degree $d$.
\end{trivlist}
\end{prop}
\par
For $d\geq 3$ let $C_{0,n,\{\star\}}(\pp,d)$ be the open substack of
${\overline
M}_{0,n+1}(\pp,d)$ consisting of maps $\mu\colon \pl\rightarrow \pp$ for
which $\mu(\pl)$ has $(d-1)(d-2)/2-1$ nodes and exactly one cusp marked
$\star$, coinciding with the $(n+1)$st marking. The strict lift of any such map
has bidegree $(d,2d-3)$. Thus
there is a lifting map
$$\lambda_C\colon C_{0,n,\{\star\}}(\pp,d)\to {\overline
M}_{0,n+1}(I,(d,2d-3)).$$
The closure of
$C_{0,n,\{\star\}}(\pp,d)$ will be denoted by $\cusplocus$; its dimension is
$3d-2+n$.
We define the {\em stack of cuspidal stable lifts} $\lcusplocus$ to be the
closure of
$\lambda_C(C_{0,n+1}(\pp,d))$ in ${\overline M}_{0,n}(\pp,(d,2d-3))$; its
dimension is likewise $3d-2+n$.
Let $e_1,\cdots,e_n$ and $e_C$ be the evaluation maps $\lcusplocus \to I$.
For each $d \geq 3$ and each $\gamma_1,\cdots,\gamma_n,\gamma_C \in A^*(I)$
we define the {\em cuspidal first-order Gromov-Witten invariant} by
$$
C_d(\gamma_1\cdots\gamma_n;\gamma_C)=\int e_1^*(\gamma_1) \cup \dots \cup
e_n^*(\gamma_n) \cup e_C^*(\gamma_C)
\cap [\lcusplocus].
$$
Suppose that among $\gamma_1,\dots,\gamma_n$ there are $a$ occurences of
the class $h^2$, also $b$ occurences of $\hd^2$, and $c$ occurences of
$h^2\hd$, where $a+b+2c=3d-1-\codim(\gamma_C)$. In this case we denote the
cuspidal Gromov-Witten invariant by $C_d(a,b,c;\gamma_C)$.
\par
\begin{thm}
\par
\begin{trivlist}
\item[(1)] $C_d(a,b,c;1)$ is the number of rational plane curves of degree
$d$, with one cusp, through $a$ general points, tangent to $b$ general
lines, and incident to $c$ general flags.
\par
\item[(2)] $C_d(a,b,c;h)$ is the number of such curves having the cusp on a
specified general line.
\par
\item[(3)] $C_d(a,b,c;h^2)$ is the number of such curves having the cusp at a
specified point.
\par
\item[(4)] $C_d(a,b,c;\hd)$ is the number of such curves for which the cusp
tangent line passes through a specified point.
\par
\item[(5)] $C_d(a,b,c;\hd^2)$ is the number of such curves for which the
cusp tangent line is a specified line.
\par
\item[(6)] $C_d(a,b,c;h^2\hd)$ is the number of such curves for which the
cusp is a specified point and the cusp tangent is a specified line through
that point.
\end{trivlist}
\end{thm}
\begin{proof}
The proof is almost identical to the proof of Theorem \ref{gwenumsign2}.
\end{proof}
\begin{prop}
\begin{trivlist}
\par
\item[(1)] For all $d\geq 3$ and all $\gamma_1,\dots,\gamma_{n-1}$,
$$
C_d(\gamma_1\cdots\gamma_{n-1}\cdot 1;\gamma_C)=0. $$
\item[(2)] If $\gamma_n$ is a divisor class in
$A^1(I)$, then for all $d\geq 3$ and all $\gamma_1,\dots,\gamma_{n-1}$,
$$
C_d(\gamma_1\cdots\gamma_n;\gamma_C)=
C_d(\gamma_1\cdots\gamma_{n-1};\gamma_C)\int (\gamma_n\cap [C]), $$
where
$[C]=d \hd^2 +(2d-3) h^2$ is the class of a the strict lift of a rational
plane curve of degree $d$ with one cusp and $(d-1)(d-2)/2-1$ nodes.
\end{trivlist}
\end{prop}
\par
%
%
%
%
\section{Special boundary divisors on the space of stable lifts}
\label{sectiondivisors}
\par
The components of the special boundary divisors of
${\overline M}_{0,n}(I,(d,2d-2))$ are indexed by the set of $4$-tuples
$\left(A_1,A_2,(d_1,c_1),(d_2,c_2)\right)$
where $d_1+d_2=d$, $c_1+c_2=2d-2$,
and $A_1\cup A_2=\{1,\dots,n\}$ is a partition in which two of the four numbers
$1,2,3,4$ belong to $A_1$ and the other two to $A_2$. The general member of the
corresponding divisor is a stable map with two components, one with
markings indexed by $A_1$ and of bidegree $(d_1,c_1)$, and the other with
markings indexed by $A_2$ and of bidegree $(d_2,c_2)$ \cite{Pandharipande}.
Since the general stable map of $\lmoduli$ is not contained in this
divisor, its restriction to $\lmoduli$ is a divisor on the space of stable
lifts, which we will denote by $D\left(A_1,A_2,(d_1,c_1),(d_2,c_2)\right)$.
This divisor may be empty; for example, $D\left(A_1,A_2,(d,0),(0,2d-2)\right)$
is empty for all $A_1$ and $A_2$.
\par
The linear equivalence of the special boundary divisors on
${\overline M}_{0,n}(I,(d,2d-2))$ passes to their restrictions:
$$
D(12\mid 34)\cong D(13\mid 24)\cong D(14\mid 23),
$$
where
\begin{equation}\label{sumpart}
D(12\mid 34)=\sum D\left(A_1,A_2,(d_1,c_1),(d_2,c_2)\right),
\end{equation}
the sum over all partitions in which $1,2\in A_1$ and
$3,4\in A_2$; the other two divisors have similar decompositions.
Thus for all choices of cohomology classes
$\gamma_1,\cdots,\gamma_n$
we have a numerical equality
\begin{equation}
\int e^*_1(\gamma_1)\cup \dots \cup e^*_n(\gamma_n)
\cap [D(12\mid 34)]
=\int e^*_1(\gamma_1)\cup \dots \cup e^*_n(\gamma_n)
\cap [D(14\mid 23)].
\end{equation}
\par
To obtain these equations in an explicit form, we must identify the
components of $D(12\mid 34)$ and $D(14\mid 23)$, at least to the extent
that they affect our calculations. We say that a divisor $D$ on $\lmoduli$
is {\em irrelevant} if for any
$\gamma_1,\dots,\gamma_n$ in $A^*(I)$ we have
\begin{equation}
\label{vanish}
\int e^*_1(\gamma_1)\cup \dots \cup e^*_n(\gamma_n)\cap[D]=0.
\end{equation}
Other divisors are said to be {\em relevant}.
\par
Let $\tau\colon\lmoduli\to {\overline
M}_{0,0}(\pp,d)$ be the map which composes stable maps $\mu\colon C\to I$
with $I\to \pp$, forgets all markings and contracts any components which
have become unstable.
\begin{prop}\label{unless}
A component $D$ of $\ldivisor$ is irrelevant unless
$\codim(\tau(D))\leq 1$ in ${\overline M}_{0,0}(\pp,d)$.
\end{prop}
\begin{proof}
It suffices to establish (\ref{vanish}) when $\gamma_1,\dots,\gamma_n$ are
elements of the basis $\{T_0,T_1,T_2,T_3,T_4,T_5\}$.
Assume first that one of the classes, say $\gamma_n=1$ is the identity.
Then by the projection formula applied to the morphism
$\pi\colon\lmoduli\to\lamoduli{n-1}{d}$ which forgets the $n$th
marking, it follows that the above degree is equal to $$
\int e^*_1(\gamma_1)\cup \dots \cup e^*_{n-1}(\gamma_{n-
1})\cap\pi_*[D]. $$
However, $\pi_*[D]$ vanishes unless
$\dim(\pi(D))=\dim(D)=\dim(\lamoduli{n-1}{d})$, so
$\pi(D)=\lamoduli{n-1}{d}$ and thus $\codim(\tau(D))=0$. \par
It remains to consider the case when none of the $\gamma_t$'s is the
identity. We will count dimensions. The dimension of $\lmoduli$ is
$3d-1+n$, and
\begin{equation*}
 \codim(e_i^*(h))=\codim(e_i^*(\hd))=1, \quad
\codim(e_i^*(h^2))=\codim(e_i^*(\hd^2))=2, \quad
\codim(e_i^*(h^2\hd))=3. \end{equation*}
Let $a$ be the number of pullbacks
of $h^2$, $b$ the number of pullbacks of $\hd^2$ and $c$ the number of
pullbacks of $h^2\hd$.
Then the class $e^*_1(\gamma_1)\cap \dots \cap
e^*_n(\gamma_n)$ has dimension 1 if  $a+b+2c=3d-2$. The number of pullbacks
of the divisor
classes $h$ and $\hd$ does not influence the vanishing nor the dimension of
the class. The effect of a divisor class on the degree of the class,
evaluated on
$[D]$ is multiplication by $d$ in the case of $h$, and multiplication by $2d-2$
in the case of $\hd$, by Proposition \ref{GWprop}. \par
Assume $\codim(\tau(D))\geq 2$; then $\dim(\tau(D))\leq 3d-3$, and a dimension
count checks that there are no stable maps in $\tau(D)$ through $a$ general
points, tangent to $b$ lines and incident to $c$ flags with $a+b+2c=3d-2$.
Hence one must put at least one condition on a map to a fiber of
$I\to \pp$ contracted by $f$. However, because of Proposition
\ref{liftclassif} this will not bring down the count $3d-2$ of conditions
on the
maps in $\tau(D)$. Indeed, putting one tangency condition on a map to a
fiber will, in case of a fiber over the image of an attachment point, induce an
extra point condition on the maps in $\tau(D)$, and in the case of a fiber over
a singularity of maps in $\tau(D)$, induce the condition that a singularity
should be on a line; also a codimension 1 condition. Similarly, putting two
tangent conditions or one flag condition on a map to a fiber will induce two
extra point conditions or the condition that a singularity of the map should be
at a fixed point, respectively. These are both codimension 2 conditions. It is
clearly impossible to put point conditions on the fiber maps, so one concludes
that the evaluation is zero unless $\codim(\tau(D))\leq 1$. \end{proof}
\begin{cor}\label{irrel}
A component $D$ of $\ldivisor$
is irrelevant unless the bidegrees $(d_1,c_1)$ and $(d_2,c_2)$ are either
\begin{trivlist}
\item[(1)] $(d_1,2d_1-2)$ and $(d_2,2d_2)$ with $d_1\geq 1$, $d_2\geq 0$
and $d_1+d_2=d$, or \item[(2)] $(d,2d-3)$ and $(0,1)$ with $d\geq 2$,
\end{trivlist}
or vice versa.
\end{cor}
\begin{proof}
First case: assume that $d_1=d$ and $d_2=0$ or vice versa.
The locus of curves in ${\overline M}_{0,0}(\pp,d)$ of degree $d$ and class
$c$ is dense if $c=2d-2$, has codimension 1 if $c=2d-3$ and codimension at
least 2 for all other values of $c$. (The case $d=2$ is special. Here $\tau(D)$
could parametrize degree two covers of lines in $\pp$.) \par
Second case: $d_1\geq 1$ and $d_2\geq 1$. We have that $\tau(D)$ is a subset of
the irreducible divisor $D(d_1,d_2)$ on ${\overline M}_{0,0}(\pp,d)$. Therefore
$D$ is irrelevant unless $\tau(D)=D(d_1,d_2)$. Now, by Proposition
\ref{stableliftd1d2} the stable lift of a generic member of $D(d_1,d_2)$ is a
map from a curve with three components and the bidegrees of the map, restricted
to the three components, are $(d_1,2d_1-2)$, $(0,2)$ and $(d_2,2d_2-2)$. Thus,
if $D$ is relevant
then the bidegrees must be $(d_1,2d_1-2)$ and $(d_2,2d_2)$, or vice versa.
\end{proof}
\par
We will next give a detailed description of the relevant divisors. To
accomplish this we need to introduce some more
special stacks of stable maps.
\par

The stack $M_{0,n+2}(I,(0,2))$ parametrizes maps of degree $2$ from $\pl$ to a
fiber of $I$ over $\pp$. Let $M_{0,n+2,\{\star,\diamond\}}(I,(0,2))$ be the
substack in which we demand that the last two markings $\star$ and $\diamond$
are the ramification points; let ${\overline
M}_{0,n,\{\star,\diamond\}}(I,(0,2))$
be its closure in ${\overline M}_{0,n+2}(I,(0,2))$.
In the same way, for maps of degree 2 to a fiber of $I$ over $\dpp$, we define
${\overline M}_{0,n,\{\star,\diamond\}}(I,(2,0))$.
\par
We define
\begin{equation}\label{emoduli}
{\overline M}^1_{0,n,\{\star\}}(\pp,(d,2d))=\bigcup_{A_1\cup A_2=\{1,\dots,n\}}
{\overline M}^1_{0,A_1\cup \{\diamond\}}(\pp,d)\times_I {\overline
M}_{0,A_2,\{\star,\diamond\}}(I,(0,2))
\end{equation}
where each fiber product is defined by using the evaluation maps corresponding
to $\diamond$; it is a closed substack of ${\overline M}_{0,n+1}(I,(d,2d))$.
The general member of each component of ${\overline
M}^1_{0,n,\{\star\}}(\pp,(d,2d))$
represents a map from a curve with two components. The restriction of the map
to one component is the strict lift of a map to a rational nodal curve; the
restriction to the other component is as described in the previous paragraph.
Similarly, we define a closed substack of ${\overline M}_{0,n+1}(I,(2,1))$ by
\begin{equation}\label{ctwomoduli}
{\overline C}^1_{0,n,\{\star\}}(\pp,2)=\bigcup_{A_1\cup A_2=\{1,\dots,n\}}
{\overline
M}_{0,A_1,\{\star,\diamond\}}(I,(2,0))\times_I {\overline
M}_{0,A_2\cup\{\diamond\}}(I,(0,1)).
\end{equation}
The general member of each component of ${\overline
C}^1_{0,n,\{\star\}}(\pp,2)$
represents a map from a curve with two components. The restriction of the map
to one component is a degree $2$ cover of a line of $\pp$ ; the restriction to
the other component is a map of degree $1$ to the pencil of all tangent
directions through the ramification point $\diamond$ of the first component.
The other ramification point is specially marked $\star$, and there are $n$
additional markings. The reason for the notation is that the
bidegree $(2,1)$ follows the pattern of bidegrees $(d,2d-3)$ for cuspidal
curves.
\par
\newpage
\begin{prop}\label{ldivisors}
Ignoring irrelevant components, we have the following isomorphisms of
divisors in $\lmoduli$:
\begin{trivlist}
\item[(1)] $D\left(A_1,A_2,(d,2d-2),(0,0)\right)\cong
{\overline M}^1_{0,A_1\cup \{\star\}}(\pp,d)
\times_I {\overline M}_{0,A_2\cup\{\star\}}(I,(0,0))$ for
$d\geq 1$.
\item[(2)] $D\left(A_1,A_2,(d,2d-3),(0,1)\right)\cong {\overline
C}^1_{0,A_1,\{\star\}}(\pp,d)\times_I{\overline
M}_{0,A_2\cup\{\star\}}(I,(0,1))$ for
$d\geq 2$.
\item[(3)] $D\left(A_1,A_2,(d_1,2d_1-2),(d_2,2d_2)\right)\cong
{\overline M}^1_{0,A_1\cup\{\star\}}(\pp,d_1)\times_I {\overline
M}^1_{0,A_2,\{\star\}}(\pp,(d_2,2d_2))$ for $d_1,d_2>0$ with $d_1+d_2=d$.
\end{trivlist}
In particular each relevant component occurs with multiplicity one.
\end{prop}
\begin{proof}
\par
(1) Note that on the right side we have a fiber product of smooth
morphisms, hence an irreducible stack.
The forgetful morphism
$$
D\left(A_1,A_2,(d,2d-2),(0,0)\right)
\to {\overline M}^1_{0,A_1\cup \{\star\}}(\pp,d)
$$
is the restriction of the projection
$$
{\overline M}_{0,A_1\cup \{\star\}}(I,(d,2d-2))
\times_I {\overline M}_{0,A_2\cup\{\star\}}(I,(0,0))
$$
onto its first factor.
Comparing dimensions, we conclude that
we have the desired isomorphism of stacks.
\par
To see that the multiplicity is one, we remark that the lifting
map $\lambda$ extends to general points of our divisor.
Indeed, a general point represents the lift of a stable map from
a curve with two components, consisting of an immersion onto a rational
degree $d$ curve together with a constant map to some nonsingular point.
The unique stable lift consists of the strict lift of the immersion
together with the constant map to the unique tangent direction at the
point. (See Figure 5.1.)
Since $\lambda$ is a birational morphism, and since
$D(A_1,A_2,d,0)$ is a reduced divisor on $\ppmoduli$,
our divisor is likewise reduced.
\par
$$
\xy
0;<0.5cm,0cm>:
(-1,-1);(1,1)**\dir{-};
(-1,1);(1,-1)**\dir{-};
(0,0)*{\bullet};
(-0.5,0.5)*{\bullet};
(0.5,-0.5)*{\bullet};
(1.5,0);(3.5,0)**\dir{-};
?>*\dir{>};
(4,0.5);(5.5,0)**\crv{(8,-0.5)&(4,-1.5)};
(5.5,0);(6,0.5)**\dir{-};
(5,0.185)*{\bullet};
(1.5,1.5);(3.5,3.5)**\dir{-};
?>*\dir{>};
(5,3.5);(5,1.5)**\dir{-};
?>*\dir{>};
(4,5.5);(5.4,4.9)**\crv{(8,4.5)&(4,3.5)};
(5.55,5.05);(6,5.5)**\dir{-};
(5,5.185)*{\bullet};
\endxy
$$
\par
\bigskip
\begin{center}
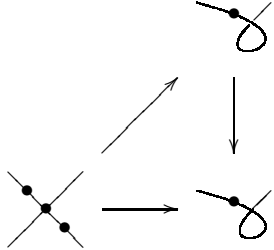
{\bf Figure 5.1} Stable lift of $D\left(A_1,A_2,(d,2d-2),(0,0)\right)$,
with two markings in $A_2$ shown.
\end{center}
\par
(2) For $d \geq 3$ we again observe that on the right side we have a fiber
product of smooth morphisms, hence  an irreducible stack. By Proposition
\ref{unless}, the image $\tau(D)$ of a relevant component $D$ of
$D\left(A_1,A_2,(d,2d-3),(0,1)\right)$ must have codimension at most $1$.
But a point of $\tau(D)$ cannot represent an immersion onto a nodal curve,
since the lift of such a curve has bidegree $(d,2d-2)$.
Nor can it represent a map to a reducible curve, since by  Proposition
\ref{stableliftd1d2} the stable lift of a general map of this type has a
component of bidegree $(0,2)$.
Hence $\tau(D)$ must be the divisor on ${\overline M}_{0,0}(\pp,d)$
representing cuspidal curves. By Proposition \ref{liftclassif}, the stable lift
of a general cuspidal curve consists of its strict lift together with a map to
the pencil of tangent directions at the cusp.
Hence we have the
desired isomorphism of stacks (ignoring irrelevant components).
\par
To see that the multiplicity is one, we first consider the case when $A_2$
is empty. By a local calculation we will verify that the lifting map
$\lambda$ extends to general points of our divisor.
The deformation theory of a node degenerating into a cusp is given in local
coordinates by $x^3+\epsilon x^2=y^2$. We may parametrize this family by
$$
(\epsilon,t)\mapsto (t^2-\epsilon,t^3-t\epsilon).
$$
Then the family of strict lifts is given by
$$
(\epsilon,t)\mapsto
\bigg((t^2-\epsilon,t^3-t\epsilon),[2t:3t^2-\epsilon]\bigg)
$$
for $\epsilon\ne 0$.
This map is not determined at $\epsilon=0$, $t=0$. Blowing up the point of
indeterminacy by introducing projective coordinates $[\epsilon_1:t_1]$
satisfying $\epsilon t_1=t\epsilon_1$, we extend the map to
$$
\bigg((\epsilon,t),[\epsilon_1:t_1]\bigg)\mapsto
\bigg((t^2-\epsilon,t^3-t\epsilon),[2t_1:3tt_1-\epsilon_1]\bigg).
$$
The special map for $\epsilon=0$ is (locally) the stable lift of $D$. It
has two components: the strict lift of the cuspidal curve, with
bidegree$(d,2d-3)$, together with a map to the pencil of tangent directions
at the cusp, with bidegree $(0,1)$.
(See Figure 5.2.)
Thus, as claimed, the lifting map $\lambda$ is defined at the point of
$D\left(A_1,\emptyset,(d,2d-3),(0,1)\right)$ representing the cuspidal
curve. And thus, by the same argument as in part (1), this divisor is
reduced.
\par
$$
\xy
0;<0.5cm,0cm>:
(-1,-1);(1,1)**\dir{-};
(-1,1);(1,-1)**\dir{-};
(0,0)*{\bullet};
(1.5,0);(3.5,0)**\dir{-};
?>*\dir{>};
(4,-1);(4,1)**\dir{-};
(6,-1);(6,1)**\crv{(2,0)};
(4,0)*{\bullet};
(0,-1.5);(0,-3.5)**\dir{-};
?>*\dir{>};
(5,-1.5);(5,-3.5)**\dir{-};
?>*\dir{>};
(-1,-6);(1,-4)**\dir{-};
(0,-5)*{\bullet};
(1.5,-5);(3.5,-5)**\dir{-};
?>*\dir{>};
(6,-4);(4,-5)**\crv{(6,-5)};
(6,-6);(4,-5)**\crv{(6,-5)};
(4,-5)*{\bullet};
\endxy
$$
\par
\bigskip
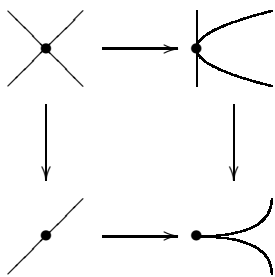
\begin{center} {\bf Figure 5.2} Stable lift of
$D\left(A_1,A_2,(d,2d-3),(0,1)\right)$, $d\geq 3$.
\end{center}
\par
Interpreting this calculation in a different way, we have shown that, for
the unique stable lift of a stable map to a cuspidal curve, we can build a
one-parameter family for which the general member is the strict lift of an
immersion and for which the special member is the specified cuspidal stable
lift; since the special member appears in the family with multiplicity one,
the divisor $D\left(A_1,\emptyset,(d,2d-3),(0,1)\right)$ must be reduced.
Similarly, we can build such a family when $A_2$ is nonempty. We assume we
are at a general point of the divisor
$D\left(A_1,A_2,(d,2d-3),(0,1)\right)$, so that the corresponding map has
just two components, with the markings indexed by $A_2$ lying on the $\pl$
mapping to the pencil of directions at the cusp. Using the same family of
maps as above, we let $L_1, L_2, \dots$ be lines through the origin in the
$\epsilon$-$t$ plane whose directions correspond to these markings, and let
the markings on nearby members of the family be the points of intersection
with these lines. Then blowing up the origin, as above, creates a family
whose special member is the specified cuspidal stable lift, with the
specified markings. Since this member appears in the family with
multiplicity one, the divisor
$D\left(A_1,A_2,(d,2d-3),(0,1)\right)$ must be reduced.
\par
We now consider the case $d=2$. The image $\tau(D)$ of a relevant component
$D$ of $D\left(A_1,A_2,(2,1),(0,1)\right)$ must have codimension at most
$1$. But a point of $\tau(D)$ cannot represent a map to a nonsingular
conic, since the lift of such a curve has bidegree $(2,2)$, nor a map to a
pair of distinct lines, since by Proposition \ref{stableliftd1d2} the
stable lift of such a map has a component of bidegree $(0,2)$. Hence
$\tau(D)$ must be the stack of degree $2$ covers of lines in $\pp$. By
Proposition \ref{liftclassif}, we have the desired isomorphism of stacks
(again ignoring irrelevant components).
\par
To see that the multiplicity is one, we consider a general point of the
divisor. This point represents the stable lift of a degree $2$ cover of a
line. A deformation of such a map is given by
$$
(\epsilon,t)\mapsto (\epsilon t,t^2),
$$
and the family of lifts is given by
$$
(\epsilon,t)\mapsto \bigg((\epsilon t,t^2),[\epsilon:2t]\bigg).$$
Blowing up the point of indeterminacy $\epsilon=t=0$ by introducing
projective coordinates $[\epsilon_1:t_1]$ satisfying $\epsilon
t_1=t\epsilon_1$, we extend the family of lifts to
 $$
\bigg((\epsilon,t),[\epsilon_1:t_1]\bigg)\mapsto \bigg((\epsilon
t,t^2),[\epsilon_1:2t_1]\bigg). $$
We do a similar blowup at the other ramification point
(located at $t = \infty$). Then the
special member $\epsilon=0$ has one component that is the strict lift of
the cover of the line and two components mapping into the pencils of
tangent directions at the ramification points.
(See Figure 5.3.)
As in the case $d \geq 3$, we can arrange for the markings on the special
member to appear in any specified configuration. Since the special member
appears in the family with multiplicity $1$, the divisor
$D\left(A_1,A_2,(2,1),(0,1)\right)$ must be reduced.
\par
$$
\xy
0;<0.5cm,0cm>:
(-1,0);(1,0)**\dir{-};
(-0.333,-1);(-0.333,1)**\dir{-};
(-0.333,0)*{\bullet};
(-0.333,0)*{\times};
(0.333,-1);(0.333,1)**\dir{-};
(0.333,0)*{\bullet};
(0.333,0)*{\times};
(1.5,0);(3.5,0)**\dir{-};
?>*\dir{>};
(4,0);(6,0)**\dir{-};
(4.667,-1);(4.667,1)**\dir{-};
(4.667,0)*{\bullet};
(4.667,0)*{\times};
(5.333,-1);(5.333,1)**\dir{-};
(5.333,0)*{\bullet};
(5.333,0)*{\times};
(0,-1.5);(0,-3.5)**\dir{-};
?>*\dir{>};
(5,-1.5);(5,-3.5)**\dir{-};
?>*\dir{>};
(-1,-5);(1,-5)**\dir{-};
(-0.333,-5)*{\times};
(0.333,-5)*{\times};
(1.5,-5);(3.5,-5)**\dir{-};
?>*\dir{>};
(4,-5);(6,-5)**\dir{-};
(4.667,-5)*{\times};
(5.333,-5)*{\times};
\endxy
$$
\par
\bigskip
\begin{center}
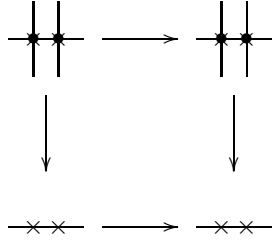
{\bf Figure 5.3}  Stable lift of $D\left(A_1,A_2,(d,2d-3),(0,1)\right)$,
$d=2$.
\end{center}
\par
(3) The image $\tau(D)$ of a relevant component $D$ of
$D\left(A_1,A_2,(d_1,2d_1-2),(d_2,2d_2)\right)$ must have codimension at
most $1$. A point of $\tau(D)$ cannot represent a map to an irreducible
curve of degree $d$. (Nor, in case $d=2$, can it represent a double cover
of a line.) Hence $\tau(D)$ must be a special boundary divisor on
${\overline M}_{0,0}(\pp,d)$:
$$
\tau(D) =
{\overline M}_{0,\{\star\}}(\pp,d_1)
\times_\pp
{\overline M}_{0,\{\star\}}(\pp,d_2)
$$
for some $d_1$ and $d_2$.
A general point of $D$ represents a
stable lift of the sort of map described in Proposition \ref{stableliftd1d2},
which tells us that the unique stable lift has three components
of bidegrees
$(d_1,2d_1-2)$, $(0,2)$ and $(d_2,2d_2-2)$, and
thus shows that  we have the desired isomorphism of stacks. (Refer again to
Figure 3.3.)
\par
To see that the divisor is reduced, we use the family of maps described in
the proof of Proposition \ref{stableliftd1d2}, with total space the surface
$xy=\epsilon^2$. The special member $\epsilon=0$ occurs with multiplicity
$1$; thus the divisor must be reduced. Furthermore we can specify markings
on the nearby members so as to obtain any configuration of specified
markings on the central component of the special member, as follows.
We introduce projective coordinates $[x_1:y_1:\epsilon_1]$, so that the
central component is the curve $x_1y_1=\epsilon_1^2$. Then to obtain
$[x_1:y_1:\epsilon_1]$ as a marking on the central component, we use
$(\frac{x_1}{\epsilon_1}\epsilon,\frac{y_1}{\epsilon_1}\epsilon,\epsilon)$
as a marking on the general member.
\end{proof}
%
%
%
%
\par
\section{Quantum potentials}\label{sectionpotentials} \par
We now introduce potential functions for each of the seven different kinds
of stable maps we are considering. We are not concerned with questions of
convergence, and all of these potentials should be interpreted as formal
power series in the indeterminates $y_0,\dots,y_5$, where
$$
\gamma=y_0T_0+\dots+y_5T_5
$$
is an arbitrary element of $A^*(I) \otimes \q$
(or in some cases as formal power series in two sets of indeterminates).
Our definitions are inspired by Proposition \ref{ldivisors} and the
auxiliary equations (\ref{emoduli}), (\ref{ctwomoduli}).
\par
We define $\calP$ to be the classical potential of the incidence
correspondence:
$$
\calP=\sum_{n \geq 3} \frac{1}{n!}
\int e^*_1(\gamma)\cup \dots \cup e^*_n(\gamma) \cap [\icmoduli].
$$
We begin the summation at $n=3$ since otherwise the moduli stacks are
empty. In fact the only nonzero term is the first one, which encodes the
intersection product:
$$
\calP=
\frac{y_0^2y_5}{2}+y_0y_1y_4+y_0y_2y_3+\frac{y_1^2y_3}{2}+\frac{y
_1y_3^2}{2}.
$$
\par
Similarly, we define a quantum potential
$$
\calN=\sum_{\substack{n \geq 0 \\ d\geq 1}}  \frac{1}{n!}
\int e^*_1(\gamma)\cup \dots \cup e^*_n(\gamma) \cap [\lmoduli].
$$
By Proposition \ref{GWprop} we have
$$
\calN=\sum_{\substack{d\geq 1 \\a +b+2c=3d-1\\ a,b,c\geq 0}}
\frac{N_d(a,b,c)y_2^ay_4^by_5^c\exp(dy_1+(2d-2)y_3)}{a!b!c!}. $$
\par
To define the quantum potential for cuspidal stable lifts we need to
introduce a second arbitrary cohomology class
$$
\delta=z_0T_0+\dots+z_5T_5.
$$
The potential is
\begin{align*}
\calC
&=\sum_{\substack{n \geq 0 \\ d\geq 2}}  \frac{1}{n!}
\int e^*_1(\gamma)\cup \dots \cup e^*_n(\gamma)
\cup e_C^*(\delta)
\cap [\lcusplocus] \\
&=\sum_{\substack{d\geq 2\\ a+b+2c=3d-2-\dim(T_s)\\a,b,c\geq
0\\s=0,\dots,5}} \frac{z_s C_d(a,b,c;T_s)y_2^ay_4^by_5^c\exp(dy_1+(2d-
3)y_3)}{a!b!c!}.
\end{align*}
\par
We will also need three potentials corresponding to maps to a
projective line.
The potential for stable maps of degree 1 to a fiber of $I$ over $\pp$ is
\begin{align*}
\calF&=\sum_{n \geq 0} \frac{1}{n!}
\int e^*_1(\gamma)\cup \dots \cup e^*_n(\gamma) \cap
[{\overline M}_{0,n}(I,(0,1))] \\
&=\left(\frac{y_4^2}{2}+y_5\right)\exp(y_3).
\end{align*}
\par
The potential for stable maps of degree 2 to such a fiber, with the
ramification at two specially marked points, is
\begin{equation*}
\calR=\sum_{n \geq 0} \frac{1}{2n!}
\int e^*_1(\gamma) \cup \dots \cup e^*_n(\gamma)
\cup e^*_{n+1}(\delta) \cup \cup e^*_{n+2}(\delta)
\cap
[{\overline M}_{0,n,\{\star,\diamond\}}(I,(0,2))].
\end{equation*}
Here it is important, in case $n=0$, that we use the {\em stack} of stable
maps of this type, since every such map has a nontrivial automorphism. To
write $\calR$ in an explicit way, we use the Gromov-Witten invariants for
such maps:
\begin{equation*}
R(a,b,c;\delta_1\cdot\delta_2)=
\int e^*_1(\gamma_1) \cup \dots \cup e^*_n(\gamma_1)
\cup e^*_{n+1}(\delta_1) \cup e^*_{n+2}(\delta_2)
\cap
[{\overline M}_{0,n,\{\star,\diamond\}}(I,(0,2))],
\end{equation*}
where $a$ of the $\gamma_t$'s equal the class $h^2$, $b$ of them equal the
class $\hd^2$, and the remaining $c$ of them equal $h^2\hd$. Unless $n=0$
we may interpret these as the number of maps satisfying $a$ point
conditions, $b$ tangency conditions and $c$ flag conditions, plus the
conditions $\delta_1$ and $\delta_2$ at the two ramification points; if
$n=0$, however, the invariant is a fraction:
\begin{align*}
R(0,2,0;T_3\cdot T_3)=&2, & R(0,0,1;T_3\cdot T_3)=&1, & R(0,1,0;T_3\cdot
T_4)=&1,\\
R(0,0,0;T_4\cdot T_4)=&\frac{1}{2}, & R(0,0,0;T_3\cdot T_5)=&\frac{1}{2}.
\end{align*}
Thus
$$
\calR=\left\{\frac{z_3^2}{2}(y_4^2+y_5)+
z_3z_4y_4+\frac{z_3z_5}{2}+\frac{z_4^2}{4}\right\}\exp(2y_3).
$$
\par
Similarly, the potential for stable maps  of degree 2 to a fiber of $I$
over the dual projective plane $\dpp$, with the ramification at two
specially marked points, is
\begin{align*}
\calL &= \sum_{n \geq 0} \frac{1}{2n!}
\int e^*_1(\gamma) \cup \dots \cup e^*_n(\gamma)
\cup e^*_{n+1}(\delta) \cup e^*_{n+2}(\delta)
\cap
[{\overline M}_{0,n,\{\star,\diamond\}}(I,(2,0))] \\
&=\left\{\frac{z_1^2}{2}(y_2^2+y_5)+ z_1z_2y_2+\frac{z_1z_5}{2}+
\frac{z_2^2}{4}\right\}\exp(2y_1).
\end{align*}
\par
Finally we have the potential $\calE$ for maps represented by the stacks
${\overline M}^1_{0,n,\{\star\}}(\pp,(d,2d))$:
\begin{align*}
\calE&=\sum_{\substack{n \geq 0 \\ d\geq 1}} \frac{1}{n!}
\int e^*_1(\gamma) \cup \dots \cup e^*_n(\gamma)
\cup e^*_{n+1}(\delta)
\cap
[{\overline M}^1_{0,n,\{\star\}}(\pp,(d,2d))] \\
&= \sum_{\substack{d\geq 1\\
a+b+2c=3d+1-\dim(T_s)\\a,b,c\geq 0\\s=0,\dots,5}}
\frac{z_sE_d(a,b,c;T_s)y_2^ay_4^by_5^c\exp(dy_1+(2d-3)y_3)}{a!b!c!},
\end{align*}
where $E_d(a,b,c;\delta)$ denotes the Gromov-Witten invariant for such maps
subject to $a$ point conditions, $b$~tangency conditions and $c$ flag
conditions, and subject to the condition
$\delta$ at the ramification point.
\par
Each of the divisors described in Proposition \ref{ldivisors} is a fiber
product
$M_1 \times_I M_2$ inside the moduli stack $\imoduli$, and each therefore
fits into a fiber square
$$
\xymatrix{
M_1 \times_I M_2  \ar[r]^\iota \ar[d]
& M_1 \times M_2 \ar[d]\\
I^{n+1}
\ar[r]^\Delta & I^{n+2}. }
$$
in which $\Delta$ is the diagonal inclusion that repeats the last factor.
The components of the vertical morphism on the right are the various
evaluation maps to $I$; the last two $e_{n+1}$ and $e_{n+2}$ become equal
when restricted to $M_1 \times_I M_2$.
By the decomposition (\ref{diag}) of the diagonal class, we have
\begin{equation}
\label{diageq}
\iota_*\left(e_1^*(\gamma_1)\cup\dots\cup e_n^*(\gamma_n)\right)
=\sum_{s=0}^5
e_1^*(\gamma_1)\cup\dots\cup e_n^*(\gamma_n)
\cup e_{n+1}^*(T_s)
\cup e_{n+2}^*(T_{5-s}).
\end{equation}
There is a partition $A_1 \cup A_2=\{1,\dots,n\}$, so that the evaluation
maps $e_t$ indexed by $t \in A_1$ factor through $M_1$ and those indexed by
$A_2$ factor through $M_2$. Taking degrees in (\ref{diageq}), we obtain the
following Proposition. (For details of this computation, see \cite[Lemma,
p.15]{Fulton2}.)
\begin{prop}\label{diagprop}
In this situation
\begin{align*}
\int e_1^*(\gamma_1)\cup\dots\cup e_n^*(\gamma_n)
& \cap [M_1 \times_I M_2] \\
& =\sum_{s=0}^5
\int
\bigcup_{t\in A_1}e_t^*(\gamma_t)
\cup e_{n+1}^*(T_s)
\cap [M_1]
\cdot
\int
\bigcup_{t\in A_2}e_t^*(\gamma_t)
\cup e_{n+2}^*(T_{5-s})
\cap [M_2].
\end{align*}
\end{prop}
\par
To use Proposition \ref{diagprop}, we need the following elementary
observation, which is a consequence of the chain rule.
\begin{prop}\label{potdiff}
Let $\gamma=y_0T_0+y_1T_1+\dots+y_5T_5$ and
$\delta=z_0T_0+z_1T_1+\dots+z_5T_5$. Suppose that $\calK$
is a power series in $y_0,\dots,y_5,z_0,\dots,z_5$ which can be written in
the following way:
$$
\calK=\sum_{\substack{n>0 \\m>0
}}\frac{K(\gamma^n;\gamma^{\prime m})}{n!m!}. $$
Let $k_1,\dots,k_r$ and $l_1,\dots,l_s$ have values in $\{0,\dots,5\}$.
Then
$$
\frac{\partial \calK}{\partial y_{k_1}\dots \partial y_{k_r}\partial
z_{l_1}\dots \partial z_{l_s}}= \sum_{\substack{n>0 \\m>0
}}\frac{K(\gamma^nT_{k_0}\dots T_{k_r};\gamma^{\prime
m}T_{l_0}\dots T_{l_r})}{n!m!}. $$
\end{prop}
\par
\begin{prop}\label{Epot}
$\displaystyle \calE=\sum_{s=0}^5
\frac{\partial \calN}{\partial y_s}
\frac{\partial \calR}{\partial z_{5-s}}.$
\end{prop}
\begin{proof}
Set
$\gamma_1=\dots=\gamma_n=\gamma=y_0T_0+y_1T_1+\dots+y_5T_5$.
Apply Proposition \ref{diagprop} to each component of
${\overline M}^1_{0,n,\{\star\}}(\pp,(d,2d))$, as listed in (\ref{emoduli}).
Sum
this identity over all components, over all $d\geq 1$, and over all $n\geq
1$. Then apply Proposition \ref{potdiff}.
\end{proof}
\par
Let $\calC_2$ be the leading term ($d=2$) of the potential $\calC$. Then a
similar argument, applied to (\ref{ctwomoduli}), yields the following
identity.
\par
\begin{prop}\label{Ctwopot}
$\displaystyle \calC_2=\sum_{s=0}^5
\frac{\partial \calL}{\partial z_s}
\frac{\partial \calF}{\partial y_{5-s}}.$
\end{prop}
\par
The same sort of argument yields the basic identity from which we will
derive recursive equations for characteristic numbers. Let $i$, $j$, $k$, and
$l$ be integers between 1 and 5. Let
\begin{equation*}
{\mathcal G}(ij\mid kl)=
\sum_{n \geq 4}  \int e_1^*(T_i) \cup e_2^*(T_j) \cup e_3^*(T_k) \cup
e_4^*(T_l) \cup e_5^*(\gamma) \cup \dots \cup e_n^*(\gamma) \cap [D(12\mid
34)].
\end{equation*}
Then, by (\ref{sumpart}), Corollary \ref{irrel} and Proposition
\ref{ldivisors},
\begin{align*}
{\mathcal G}(ij\mid kl)=
\sum_{s=0}^5 \bigg\{ &\frac{\partial^3
\calN}{\partial y_i \partial y_j \partial y_s}\frac{\partial^3 \calP}{\partial
y_k\partial y_l \partial y_{5-s}}
+\frac{\partial^3 \calN}{\partial y_k
\partial y_l \partial y_s}\frac{\partial^3 \calP}{\partial y_i\partial y_j
\partial y_{5-s}}\\
+&\frac{\partial^3 \calC}{\partial y_i \partial y_j \partial
z_s}\frac{\partial^3 \calF}{\partial y_k\partial y_l \partial y_{5-s}}
+\frac{\partial^3 \calC}{\partial y_k \partial y_l \partial
z_s}\frac{\partial^3 \calF}{\partial y_i\partial y_j \partial y_{5-s}}\\
+&\frac{\partial^3 \calN}{\partial y_i \partial y_j
\partial y_s}\frac{\partial^3 \calE}{\partial y_k\partial y_l \partial z_{5-s}}
+\frac{\partial^3 \calN}{\partial y_k \partial y_l \partial
y_s}\frac{\partial^3 \calE}{\partial y_i\partial y_j \partial z_{5-s}}
\bigg\}.
\end{align*}
The linear equivalence of the divisors $D(12\mid 34)$
and $D(13\mid 24)$ immediately implies the basic identity,
which we record as a Theorem.
\par
\begin{thm}\label{ijklrel}
For each $i,j,k$ and $l$ in $\{1,2,3,4,5\}$, there is an identity $$
\calG(ij \mid kl)=\calG(il \mid jk).
$$
\end{thm}
\par
%
%
%
%
\section{Recursive formulas for characteristic
numbers}\label{sectionrecursion}
\par
In this section we will
derive recursive
formulas for the characteristic numbers $N_d(a,b,c)$,
$C_d(a,b,c;1)$, $C_d(a,b,c;h)$, and $C_d(a,b,c;h^2)$.
For characteristic numbers of the first type, the base cases are
$$
N_1(2,0,0)=N_1(0,0,1)=1 \quad \text{and} \quad N_1(0,2,0)=0.
$$
\par
The cuspidal characteristic numbers make sense only for $d \geq 2$.
In the case $d=2$  they are defined using the stack
${\overline C}^1_{0,n,1}(\pp,2)$ whose decomposition is shown in
(\ref{ctwomoduli}).
Proposition \ref{Ctwopot} tells us how to calculate the corresponding
potential $\calC_2$ from the known potentials $\calL$ and $\calF$.
Equating coefficients, we find that the only nonzero cuspidal
characteristic numbers in degree $d=2$ are
\begin{align*}
C_2(2,1,0;T_1)=&2, & C_2(0,1,1;T_1)=&1, & C_2(1,2,0;T_1)=&1, &
C_2(1,0,1;T_1)=1,\\ C_2(1,1,0;T_2)=&1, & C_2(0,2,0;T_2)=&\frac{1}{2}, &
C_2(0,0,1;T_2)=&\frac{1}{2}, & C_2(0,1,0;T_5)=\frac{1}{2}.
\end{align*}
The fact that some of these characteristic numbers are fractions reflects
the fact that the generic map of this type has a nontrivial automorphism.
One could also calculate these characteristic numbers directly, again
taking into account the automorphisms.
\par
In addition to the desired characteristic numbers, we will also calculate
the characteristic numbers $E_d(a,b,c;T_s)$ for maps represented by the
stacks ${\overline M}^1_{0,n,1}(\pp,(d,2d))$, which appear as coefficients
in the potential $\calE$. We are not particularly interested in these
numbers, but using them simplifies many of the formulas.
\par
\begin{prop}\label{ENrels}
\begin{align*}
&E_d(a,b,c;\hd)=dbN_d(a,b-1,c)+b(b-1)N_d(a+1,b-2,c)+cN_d(a+1,b,c-1)
\text{ for $a+b+2c=3d$;}\\
&E_d(a,b,c;\hd^2)=\frac{d}{2}N_d(a,b,c)+bN_d(a+1,b-1,c)
\text{ for $a+b+2c=3d-1$;}\\
&E_d(a,b,c;h^2\hd)=\frac{1}{2}N_d(a+1,b,c)\text{ for $a+b+2c=3d-2$;}\\
&E_d(a,b,c;T_s)=0 \text{ if $s = 0$, $1$ or $2$.}
\end{align*}
\end{prop}
\begin{proof}
Use the identity of Proposition \ref{Epot}; equate the coefficients.
\end{proof}
\par
\begin{thm} \label{algthm} There is a recursive algorithm, based on Theorem
\ref{ijklrel}, for
calculating all the characteristic numbers $N_d(a,b,c)$,
$C_d(a,b,c;1)$, $C_d(a,b,c;h)$, $C_d(a,b,c;h^2)$,
$E_d(a,b,c;\hd)$, $E_d(a,b,c;\hd^2)$, and $E_d(a,b,c;h^2\hd)$.
 \end{thm}
Pandharipande \cite[Proposition 4]{Pandharipande2} proved that there is an
explicit algorithm for calculating all the numbers $N_d(a,b,0)$ and gave
explicit formulas for the numbers $C_d(3d-2,0,0;1)$ \cite[Proposition
5]{Pandharipande2}. \par
\begin{proof}
The proof will use induction on $d$.
Assume that all characteristic numbers have been
determined for degrees less than $d$. Note that the last three types of
characteristic numbers (the $E_d$'s) are determined by the other types;
hence we need only to worry about the $N_d$'s and $C_d$'s.
To simplify the argument we will denote
by $[d-1]$ any expression involving characteristic numbers for curves of
degree $d-1$ or less.
\par
First we list various cases of Theorem \ref{ijklrel} which immediately
determine a
characteristic number in degree $d$; each of these identities can be
written in the form
$$
N_d(a,b,c)=[d-1].
$$
\par
Equation $1122$ (meaning this case of the identity in Theorem \ref{ijklrel})
determines $N_d(a,b,c)$ for $a\geq 3$.
\par
Equation $1224$ determines $N_d(a,b,c)$ for $a\geq 2$ and $c\geq 1$.
\par
Equation $1155$ determines $N_d(a,b,c)$ for $a\geq 1$ and $c\geq 2$.
\par
Equation $1455$ determines $N_d(a,b,c)$ for $c\geq 3$.
\par
This leaves only six of the numbers $N_d(a,b,c)$ undetermined, namely
\begin{align*}
&N_d(2,3d-3,0), & &N_d(1,3d-4,1), & &N_d(0,3d-5,2),\\ &N_d(1,3d-2,0), &
&N_d(0,3d-3,1), & &N_d(0,3d-1,0). \end{align*}
\par
Similarily (for $d \geq 3$) there are equations of the form: $$
C_d(a,b,c;h^2)=[d-1].
$$
\par
Equation $2445$ determines $C_d(a,b,c;h^2)$ for $a\geq 1$ and $c\geq
1$.
\par
Equation $4455$ determines $C_d(a,b,c;h^2)$ for $a\geq 2$ and $c\geq
1$.
\par
Equation $2244$ determines $C_d(a,b,c;h^2)$ for $a\geq 2$.
\par
This leaves
\begin{align*}
&C_d(1,3d-5,0;h^2), & &C_d(0,3d-6,1;h^2), & & C_d(0,3d-4,0;h^2) \end{align*}
undetermined.
\par
Next, we exhibit in matrix form 8 independent equations involving the
undetermined characteristic numbers listed above, except $N_d(0,3d-1,0)$.
$$
\left(\begin{smallmatrix} 2(3d-3) & -d & 0 & 0 & 0 & 0 & 0 & 0\\
0 & 3d-4 & -2d & 0 & 0 & 0 & 0 & 0\\
2(d-2) & 0 & 0 & -d(3d-2) & d^2 & 0 & 0 & 0\\ 0 & 0 & 2(3d-5) & 0 & 0 & 0 & -d
& 0\\
0 & 0 & 2(3d-5) & 0 & 0 & -(3d-5) & 0 & 0\\ 0 & 3d-3 & -4d+4 & 0 & (3d-4)(3d-3)
& 0 & 0 & -(3d-4)d\\ 0 & d-2 & 0 & 0 & -d(3d-3) & 0 & 0 & d^2\\ 0 & 0 & 0 & 0
& -(3d-4)d & d^2 & 0 & 0\\ \end{smallmatrix}\right) \left(\begin{smallmatrix}
N_d(2,3d-3,0)\\ N_d(1,3d-4,1)\\ N_d(0,3d-5,2)\\ N_d(1,3d-2,0)\\ N_d(0,3d-
3,1)\\ C_d(1,3d-5,0;h^2)\\ C_d(0,3d-6,1;h^2)\\ C_d(0,3d-4,0;h^2)
\end{smallmatrix}
\right) = [d-1]
$$
The rows of this equation correspond to the following equations:
1124, 1145, 1123, 1445, 2345, 1345, 1135, and 1135 again (for a different
choice of $a$, $b$, and $c$).
The determinant of the matrix is
$$
- 12 d^6(3 d - 5) (d - 1) (3 d - 4) (3 d - 2),
$$
which is nonzero for all integers
$d\geq 2$. Thus these equations determine the eight degree $d$
characteristic numbers listed in the column vector.
\par
Next, the equation $1134$ can be written
\begin{align*}
C_d(a,b,c;h)&=-C_d(a,b-1,c;h^2)b
+\frac{1}{d^2}\bigg\{N_d(a,b,c+1)d(2d-2)\\
&+N_d(a+1,b+1,c)(2-d) +N_d(a,b+2,c)d\bigg\}+[d-1].
\end{align*}
It determines all the numbers $C_d(a,b,c;h)$ except $C_d(0,3d-3,0;h)$.
\par
Equation $1133$ can be written
\begin{align*}
C_d(a,&b,c;1)=-C_d(a,b-1,c;h)b-C_d(a,b,c-1;h^2)c-C_d(a,b-
2,c;h^2)\binom{b}{2}\\ &+\frac{1}{d^2}\bigg\{4N_d(a+1,b,c)(d-
1)+N_d(a,b+1,c)(3d^2-4d)\bigg\}+[d-1]\end{align*}
It determines all the numbers $C_d(a,b,c;1)$ except $C_d(0,3d-2,0;1)$.
\par
It remains to determine the three numbers $N_d(0,3d-1,0)$, $C_d(0,3d-3,0;h)$
and $C_d(0,3d-2,0;1)$.
The two equations $1133$ and $1134$ just listed give two equations. Another
is given by $3344$. To show that they are linearly independent we write the
equations in matrix form and compute the determinant.
\begin{equation*}
\left(\matrix (3d-1)d(4-3d) & d^2 & 1\\
-(3d-2)(3d-1)d & d^2 & 0\\
(3d-3)(3d-2)(3d-1) & (3d-3)(7d+2) & (3d-3)(3d-2) \endmatrix\right)
\left(\matrix N_d(0,3d-1,0)\\ C_d(0,3d-3,0;h)\\ C_d(0,3d-2,0;1)
\endmatrix \right)=
\left(\matrix \text{known} \\ \text{known} \\ \text{known} \endmatrix
\right)
\end{equation*}
Here, ``known" means expressions involving already determined numbers of
degree $d$ or less. The determinant of the matrix is
$$
6 d (d - 1) (3 d - 1) (3 d - 2) (d^2 - 4 d - 1), $$
which is nonzero for all integers $d\geq 2$. This shows that we can find the
characteristic numbers in degree $d$ from those in characteristic $d-1$.
\end{proof}
We will next give an explicit recursive algorithm, developed and
implemented by the first author, which determines all the
characteristic numbers of Theorem \ref{algthm}. The algorithm will use 7
equations compared to the 17 used in the proof of Theorem \ref{algthm}.
A drawback is, however, that the recursion is not as transparent as the
recursion of the proof. Some degree $d+1$ characteristic numbers must be
computed before all numbers in degree $d$ may be computed. \par
Using equation $1122$ we derive
the following formula, which we call equation $1122a$. We must assume that
$a\geq3$. In all sums
$d_1,d_2>0$ and $a_1,a_2,b_1,b_2,c_1,c_2\geq 0$. \begin{align*}
&N_d(a,b,c)=\\
&\sum_{\substack{d_1+d_2=d\\ a_1+a_2=a-
1\\b_1+b_2=b\\c_1+c_2=c}}
N_{d_1}(a_1,b_1,c_1)E_{d_2}(a_2,b_2,c_2;\hd^2) \left[2d_1^2d_2\binom{a-
3}{a_1-1}
-d_1^3d_2\binom{a-3}{a_1}
-d_1d_2^2\binom{a-3}{a_1-1}
\right]\binom b{b_1}\binom c{c_1}\\
&+\sum_{\substack{d_1+d_2=d\\ a_1+a_2=a\\b_1+b_2=b\\c_1+c_2=c}}
N_{d_1}(a_1,b_1,c_1)E_{d_2}(a_2,b_2,c_2;\hd) \left[2d_1d_2\binom{a-
3}{a_1-2}
-d_1^2\binom{a-3}{a_1-1}
-d_2^2\binom{a-3}{a_1-2}
\right]\binom b{b_1}\binom c{c_1}\\
\end{align*}
If we put $a=3d-1$ and $b=c=0$ we recover Kontsevich's formula \cite[Claim
5.2.1]{KontsevichManin}. If we put $c=0$ we recover the recursive formula of
Francesco and C. Itzykson \cite[2.95, p.104]{FrancescoItzykson}.
\par
The next equation $1122b$ is a rewrite of $1122a$ in a very specific case.
We use
$1122a$ with $d$ replaced by $d+1$, $a=3$, $b=3d-1$ and $c=0$. Then we solve
for $N_d(0,b,0)$, which can only occur in the first sum of the right hand side
of $1122a$. There are two cases:
\begin{align*}
&(d_1,a_1,b_1,c_1)=(d,0,b,0), &(d_2,a_2,b_2,c_2)=(1,2,0,0), \end{align*}
and the same values with the indexes 1 and 2 interchanged. Using
Proposition \ref{ENrels} we find
\begin{equation*}
E_1(2,0,0;\hd^2)=\frac{1}{2},
\quad
E_d(0,b,0;\hd^2)=\frac{d}{2}N_d(0,b,0)+bN_d(1,b-1,0).
\end{equation*}
Then we have the following equation called 1122b:
\begin{align*}
&N_d(0,b,0)=\frac{1}{d^3}\bigg\{-N_{d+1}(3,b,0)-d^2bN_d(1,b-1,0)\\
&+\sum_{\substack{d_1+d_2=d+1\\ a_1+a_2=2\\b_1+b_2=b\\
(d_i,a_i,b_i)\neq (d,0,b)}}
N_{d_1}(a_1,b_1,0)E_{d_2}(a_2,b_2,0;\hd^2) \left[2d_1^2d_2\binom{0}{a_1-
1}
-d_1^3d_2\binom{0}{a_1}
-d_1d_2^2\binom{0}{a_1-1}
\right]\binom b{b_1}\\
&+\sum_{\substack{d_1+d_2=d+1\\ a_1+a_2=3\\b_1+b_2=b}}
N_{d_1}(a_1,b_1,c_1)E_{d_2}(a_2,b_2,c_2;\hd) \left[2d_1d_2\binom{0}{a_1-
2}
-d_1^2\binom{0}{a_1-1}
-d_2^2\binom{0}{a_1-2}
\right]\binom b{b_1}\\
\end{align*}
\par
Equation $1155$: If $a\geq1$ and $c\geq 2$ then \begin{align*}
&N_d(a,b,c)=\\
&\sum_{\substack{d_1+d_2=d\\ a_1+a_2=a-1\\b_1+b_2=b\\c_1+c_2=c}}
N_{d_1}(a_1,b_1,c_1)E_{d_2}(a_2,b_2,c_2;\hd^2) \left[2d_1^2d_2\binom{c-
2}{c_1-1}
-2d_1^2d_2\binom{c-2}{c_1-2}
-d_1^3\binom{c-2}{c_1}
\right]\binom {a-1}{a_1}\binom b{b_1}\\
&+\sum_{\substack{d_1+d_2=d\\ a_1+a_2=a\\b_1+b_2=b\\c_1+c_2=c}}
N_{d_1}(a_1,b_1,c_1)E_{d_2}(a_2,b_2,c_2;\hd) \bigg[2d_1d_2\binom {c-
2}{c_1-1}
-d_1^2\binom {c-2}{c_1}
-d_2^2\binom {c-2}{c_1-2}\bigg]\binom{a-1}{a_1-1}\binom b{b_1}
\end{align*}
\par
Equation $1123a$:
If $c\geq 1$ then
\begin{align*}
N_d(a,b&,c)=\frac{1}{d^2}\bigg\{N_d(a+1,b+1,c-1)-N_d(a+2,b,c-1)(d-2)\\
+&\sum_{\substack{d_1+d_2=d\\
a_1+a_2=a+1\\b_1+b_2=b\\c_1+c_2=c-1}}
N_{d_1}(a_1,b_1,c_1)E_{d_2}(a_2,b_2,c_2;\hd^2) \left[2d_1d_2^2\binom
{a}{a_1-1}
-2d_1^2d_2\binom {a}{a_1}\right]\binom{b}{b_1}\binom{c-1}{c_1}\\
&\sum_{\substack{d_1+d_2=d\\ a_1+a_2=a+2\\b_1+b_2=b\\c_1+c_2=c-
1}} N_{d_1}(a_1,b_1,c_1)E_{d_2}(a_2,b_2,c_2;\hd)
\bigg[2d_2^2\binom{a}{a_1-2}-2d_1d_2\binom{a}{a_1-1}\bigg]
\binom{b}{b_1}\binom{c-1}{c_1} \bigg\}
\end{align*}
\par
Equation $1123b$: This is just a rewrite of 1123a. If $a\geq 1$ and $b\geq 1$
then \begin{align*}
N_d(a,b&,c)=(d-2)N_d(a+1,b-1,c)+d^2N_d(a-1,b-1,c+1)\\
+&\sum_{\substack{d_1+d_2=d\\ a_1+a_2=a\\b_1+b_2=b-
1\\c_1+c_2=c}} N_{d_1}(a_1,b_1,c_1)E_{d_2}(a_2,b_2,c_2;\hd^2)
\left[2d_1^2d_2\binom {a-1}{a_1}-2d_1d_2^2\binom {a-1}{a_1-1}
\right]\binom{b-1}{b_1}\binom c{c_1}\\
&\sum_{\substack{d_1+d_2=d\\ a_1+a_2=a+1\\b_1+b_2=b-
1\\c_1+c_2=c}} N_{d_1}(a_1,b_1,c_1)E_{d_2}(a_2,b_2,c_2;\hd)
\bigg[2d_1d_2\binom{a-1}{a_1-1}-2d_2^2\binom{a-1}{a_1-2}\bigg]
\binom{b-1}{b_1}\binom c{c_1} \end{align*}
\par
Equation $2245$: This is not an obvious recursion equation.
Given $a\geq 2,b\geq 1,c\geq 1$ and $d$ with $a+b+2c=3d-2$, we have that
\begin{align*}
0&=\sum_{\substack{d_1+d_2=d\\ a_1+a_2=a\\b_1+b_2=b\\c_1+c_2=c}}
N_{d_1}(a_1,b_1,c_1)E_{d_2}(a_2,b_2,c_2;\hd^2)d_1 \bigg[\binom{a-2}{a_1-
2}\binom {b-1}{b_1}\binom{c-1}{c_1} +\binom{a-2}{a_1}\binom {b-1}{b_1-
1}\binom{c-1}{c_1-1}\\ &-\binom{a-2}{a_1-1}\binom {b-1}{b_1}\binom{c-
1}{c_1-1} -\binom{a-2}{a_1-1}\binom {b-1}{b_1-1}\binom{c-1}{c_1}\bigg]\\
&+\sum_{\substack{d_1+d_2=d\\
a_1+a_2=a+1\\b_1+b_2=b\\c_1+c_2=c}}
N_{d_1}(a_1,b_1,c_1)E_{d_2}(a_2,b_2,c_2;\hd) \bigg[\binom{a-2}{a_1-
3}\binom {b-1}{b_1}\binom{c-1}{c_1} +\binom{a-2}{a_1-1}\binom{b-1}{b_1-
1}\binom{c-1}{c_1-1}\\ &-\binom{a-2}{a_1-2}\binom {b-1}{b_1}\binom{c-
1}{c_1} -\binom{a-2}{a_1}\binom{b-1}{b_1-1}\binom{c-1}{c_1}\bigg].\\
\end{align*}
To get a recursion equation, we rewrite Equation 2245 with $d$ replaced by
$d+1$, $a=2$
and $b+c=3d-1$. Then we solve for $N_d(0,b,c)$, which can only occur in the
first sum of the right hand side of $2245$. There are two cases:
\begin{align*}
&(d_1,a_1,b_1,c_1)=(d,0,b,c), &(d_2,a_2,b_2,c_2)=(1,2,0,0), \end{align*}
and the same values with the indexes 1 and 2 interchanged. Using
Proposition \ref{ENrels} we find
\begin{equation*}
E_1(2,0,0;\hd^2)=\frac{1}{2},
\quad
E_d(0,b,c;\hd^2)=\frac{d}{2}N_d(0,b,c)+bN_d(1,b-1,c).
\end{equation*}
We have, for $b\geq 1$, $c\geq 1$, equation 2245:
\begin{align*}
&N_d(0,b,c)=\frac{1}{d}\bigg\{-bN_d(1,b-1,c)-\\
&\sum_{\substack{d_1+d_2=d+1\\
a_1+a_2=2\\b_1+b_2=b\\c_1+c_2=c\\(d_i,a_i,b_i,c_i)\neq (d,0,b,c)}}
N_{d_1}(a_1,b_1,c_1)E_{d_2}(a_2,b_2,c_2;\hd^2)d_1 \bigg[\binom{0}{a_1-
2}\binom {b-1}{b_1}\binom{c-1}{c_1} +\binom{0}{a_1}\binom {b-1}{b_1-
1}\binom{c-1}{c_1-1}\\ &-\binom{0}{a_1-1}\binom {b-1}{b_1}\binom{c-
1}{c_1-1} -\binom{0}{a_1-1}\binom {b-1}{b_1-1}\binom{c-1}{c_1}\bigg]\\
&+\sum_{\substack{d_1+d_2=d+1\\
a_1+a_2=3\\b_1+b_2=b\\c_1+c_2=c}}
N_{d_1}(a_1,b_1,c_1)E_{d_2}(a_2,b_2,c_2;\hd) \bigg[\binom{0}{a_1-
3}\binom {b-1}{b_1}\binom{c-1}{c_1} +\binom{0}{a_1-1}\binom{b-1}{b_1-
1}\binom{c-1}{c_1-1}\\ &-\binom{0}{a_1-2}\binom{b-1}{b_1}\binom{c-
1}{c_1} -\binom{0}{a_1}\binom{b-1}{b_1-1}\binom{c-1}{c_1}\bigg] \bigg\}
\end{align*}
\par
Now we are ready to present the algorithm for calculating the numbers
$N_d(a,b,c)$ with $a+b+2c=3d-1$. The inputs to the algorithm consist of the
trivial degree 1 numbers, which are zero except for $N_1(2,0,0)=1$ and
$N_1(0,0,1)=1$. It is also understood that the instruction ``use equation
$ijkl$'' means to use this equation and then to replace all occurences of
the $E_d$'s by means of Proposition \ref{ENrels}.
\par
\begin{tabbing}
$N_d(a,b,c):=$ \= {\bf if} $(d=1)$ {\bf then} \= {\bf if} $(a=2$ {\bf and $b=0$
and $c=0)$ then $1$} \\ \> \> {\bf else if $(a=0$ and $b=0$ and $c=1)$ then
$1$ }\\ \> \> {\bf else} 0 \\
\> {\bf else if $\big((a\geq 4)$ or $(a=3$ and $c\geq 1)\big)$ then} use
equation 1122a\\
\> {\bf else if $(a\geq 1$ and $c\geq 2)$ then} use equation 1155\\ \> {\bf
else if $(c\geq 3)$ then} use equation 1123a\\ \> {\bf else if $\big((a=3$ and
$c=0)$ or $(a=2$ and $c=1)$ or $(a=1$ and $c=1)$}\\ \> \> {\bf or $(a=2$ and
$c=0)$ or $(a=1$ and $c=0)\big)$ then} use equation 1123b\\
\> {\bf else if $\big((a=0$ and $c=2)$ or $(a=0$ and $c=1)\big)$ then} use
equation 2245\\ \> {\bf else} use equation 1122b.
\end{tabbing}
\par
We finally present tables of characteristic numbers for $d \leq 5$.
In the tables of characteristic numbers $N_d(a,b,c)$, columns have the same
$a$, rows have the same $c$, and $b$ is given by
$3d-1-a-2c$.
\begin{center}
\par
\begin{tabular}{|r|c|c|c|c|c|c|} \hline
$c\downarrow$ & \multicolumn{6}{|c|}{$N_2(a,b,c)$} \\ \hline 2 & 1 & 1
&\multicolumn{4}{c|}{} \\ \cline{1-5} 1 & 1 & 2 & 2 & 1 &
\multicolumn{2}{c|}{} \\ \hline 0 & 1 & 2 & 4 & 4 & 2 & 1 \\ \hline
$a\to$ & 0 & 1 & 2 & 3 & 4 & 5 \\ \hline \end{tabular}
\par
\smallskip
\par
{\bf Table 7.1} Characteristic numbers for plane conics.
\end{center}
\par
\bigskip
\begin{center}
\begin{tabular}{|r|c|c|c|c|c|c|c|c|c|} \hline $c\downarrow$ &
\multicolumn{9}{|c|}{$N_3(a,b,c)$} \\ \hline 4 & 4 & \multicolumn{8}{c|}{}
\\ \cline{1-4} 3 & 16 & 12 & 6 & \multicolumn{6}{c|}{} \\ \cline{1-6} 2 & 56
& 56 & 40 & 20 & 8 & \multicolumn{4}{c|}{} \\ \cline{1-8} 1 & 148 & 200 &
196 & 136 & 68 & 28 & 10 & \multicolumn{2}{c|}{} \\ \hline 0 & 400 & 600 &
756 & 712 & 480 & 240 & 100 & 36 & 12 \\ \hline $a\to$ & 0 & 1 & 2 & 3 & 4
& 5 & 6 & 7 & 8 \\ \hline \end{tabular}
\par
\smallskip
\par
{\bf Table 7.2} Characteristic numbers for plane rational nodal cubics.
\end{center}
\bigskip
\par
The numbers for nodal cubics were calculated by Zeuthen \cite{Zeuthen},
Maillard \cite{Maillard}, Schubert \cite{Schubert}, Sacchiero \cite{Sacchiero},
Kleiman and Speiser \cite{KleimanSpeiser}, Aluffi \cite{Aluffi2} and
Pandharipande \cite{Pandharipande2}.
\begin{center}
\par
\begin{tabular}{|r|c|c|c|c|c|c|c|c|c|c|c|c|} \hline $c\downarrow$ &
\multicolumn{12}{|c|}{$N_4(a,b,c)$} \\ \hline 5 & 120 & 60 &
\multicolumn{10}{c|}{} \\ \cline{1-5} 4 & 816 & 528 & 264 & 108 &
\multicolumn{8}{c|}{} \\ \cline{1-7} 3 & 5040 & 3960 & 2472 & 1224 & 504 &
180 & \multicolumn{6}{c|}{} \\ \cline{1-9} 2 & 26408 & 25352 & 19424 &
11840 & 5816 & 2408 & 872 & 284 & \multicolumn{4}{c|}{} \\ \cline{1-11}
1 & 124592 & 140912 & 130824 & 97496 & 58208 & 28392 & 11792 & 4304 & 1416
& 428 & \multicolumn{2}{c|}{} \\ \hline 0 & 581904 & 728160 & 783584 &
699216 & 505320 & 295544 & 143040 & 59424 & 21776 & 7200 & 2184 & 620 \\
\hline $a\to$ & 0 & 1 & 2 & 3 & 4 & 5 & 6 & 7 & 8 & 9 & 10 & 11 \\ \hline
\end{tabular}
\par
\smallskip
{\bf Table 7.3} Characteristic numbers for plane rational nodal quartics.
\end{center}
\par
\bigskip
The numbers for nodal quartics with $c=0$ were calculated by Pandharipande
\cite{Pandharipande2}.
\par
\begin{center}
\par
\begin{tabular}{|r|c|c|c|c|c|c|c|} \hline $c\downarrow$ &
\multicolumn{7}{|c|}{$N_5(a,b,c)$} \\ \hline 7 & 840 &
\multicolumn{6}{c|}{} \\ \cline{1-4} 6 & 9120 & 4560 & 1920 &
\multicolumn{4}{c|}{} \\ \cline{1-6} 5 & 88560 & 52800 & 26160 & 11040 &
4080 & \multicolumn{2}{c|}{} \\ \hline 4 & 792432 & 548064 & 318432 &
155904 & 65712 & 24432 & 8184 \\ \hline 3 & 6347808 & 5092128 & 3442536 &
1961952 & 950880 & 400128 & 149448 \\ \hline
2 & 46200256 & 42546112 & 33296896 & 22024720 & 12343168 & 5928736 &
2489872 \\ \hline
1 & 317706976 & 327704960 & 292474400 & 222844672 & 144185600 & 79521536
& 37862920 \\ \hline
0 & 2150306368 & 2414524160 & 2397491872 & 2069215552 & 1532471744 &
969325888 & 526105120 \\ \hline
$a\to$ & 0 & 1 & 2 & 3 & 4 & 5 & 6 \\ \hline \end{tabular}
\par
\smallskip
\begin{tabular}{|r|c|c|c|c|c|c|c|c|} \hline $c\downarrow$ &
\multicolumn{8}{|c|}{$N_5(a,b,c)$} \\ \hline 3 & 50496 & 15672
& \multicolumn{6}{c|}{} \\ \cline{1-5} 2 & 932656 & 316960 & 99088 & 28816
& \multicolumn{4}{c|}{} \\ \cline{1-7} 1 & 15857120 & 5946992 & 2028160 &
636896 & 186080 & 51040 &
\multicolumn{2}{c|}{} \\ \hline
0 & 248204432 & 103544272 & 38816224 & 13258208 & 4173280 & 1222192 &
335792 & 87304 \\ \hline $a\to$ & 7 & 8 & 9 & 10 & 11 & 12 & 13 & 14 \\
\hline \end{tabular}
\par
\smallskip
{\bf Table 7.4} Characteristic numbers for plane rational nodal quintics.
\end{center}
\par
\bigskip
The numbers with $c=0$ and $a\geq 3$ were calculated by Francesco and
Itzykson \cite[2.97, p.104]{FrancescoItzykson}. \par
We use the same kind of equations to
derive the characteristic numbers $C_d(a,b,c;1)$, $C_d(a,b,c;h)$ and
$C_d(a,b,c;h^2)$.
\par
Equation $1144$:
For any $a$, $b$ and $c$ with $a+b+2c=3d-4$ we have that
\begin{align*}
&C_d(a,b,c;h^2)=\frac{1}{d^2}\bigg\{2dN_d(a,b+1,c+1)-N_d(a+1,b+2,c)\\
&+\sum_{\substack{d_1+d_2=d\\
a_1+a_2=a\\b_1+b_2=b+2\\c_1+c_2=c}}
N_{d_1}(a_1,b_1,c_1)E_{d_2}(a_2,b_2,c_2;\hd^2)
\left[2d_1^2d_2\binom{b}{b_1-1}-d_1d_2^2\binom{b}{b_1-2}
-d_1^3\binom{b}{b_1}\right]\binom {a}{a_1}\binom c{c_1}\\
&+\sum_{\substack{d_1+d_2=d\\
a_1+a_2=a+1\\b_1+b_2=b+2\\c_1+c_2=c}}
N_{d_1}(a_1,b_1,c_1)E_{d_2}(a_2,b_2,c_2;\hd)
\bigg[2d_1d_2\binom{b}{b_1-1}
-d_1^2\binom{b}{b_1}-d_2^2\binom{b}{b_1-2} \bigg]\binom{a}{a_1-
1}\binom c{c_1}\bigg\} \end{align*}
\par
In the tables of the numbers $C_d(a,b,c;h^2)$ the number of tangent conditions
$b$ is given by $3d-4-a-2c$.
\par
\begin{center}
\par
\begin{tabular}{|r|c|c|c|c|c|c|} \hline
$c\downarrow$ & \multicolumn{6}{|c|}{$C_3(a,b,c;h^2)$} \\ \hline 2 & 4 &
2 & \multicolumn{4}{c|}{} \\ \cline{1-5} 1 & 14 & 12 & 6 & 2 &
\multicolumn{2}{c|}{} \\ \hline 0 & 32 & 44 & 38 & 20 & 8 & 2 \\ \hline
$a\to$ & 0 & 1 & 2 & 3 & 4 & 5 \\ \hline \end{tabular}
\par
\smallskip
{\bf Table 7.5} Characteristic numbers for plane rational 1-cuspidal cubics;
cusp at a specified point.
\end{center}
\par
\bigskip
The numbers for cubics were calculated by Schubert \cite{Schubert}, by Kleiman
and Speiser \cite{KleimanSpeiser2}, and by Aluffi \cite{Aluffi2}.
\par
\begin{center}
\par
\begin{tabular}{|r|c|c|c|c|c|c|c|c|c|} \hline $c\downarrow$ &
\multicolumn{9}{|c|}{$C_4(a,b,c;h^2)$} \\ \hline 4 & 18 &
\multicolumn{8}{c|}{} \\ \cline{1-4} 3 & 132 & 72 & 30 &
\multicolumn{6}{c|}{} \\ \cline{1-6} 2 & 816 & 580 & 312 & 132 & 46 &
\multicolumn{4}{c|}{} \\ \cline{1-8} 1 & 4084 & 3760 & 2636 & 1420 & 616 &
224 & 70 & \multicolumn{2}{c|}{} \\ \hline 0 & 17444 & 19912 & 17904 &
12392 & 6700 & 2964 & 1112 & 364 & 102 \\ \hline $a\to$ & 0 & 1 & 2 & 3 & 4
& 5 & 6 & 7 & 8 \\ \hline \end{tabular}
\par
\smallskip
{\bf Table 7.6} Characteristic numbers for plane rational 1-cuspidal quartics;
cusp at a specified point.
\end{center}
\par
\bigskip
\begin{center}
\par
\begin{tabular}{|r|c|c|c|c|c|c|} \hline
$c\downarrow$ & \multicolumn{6}{|c|}{$C_5(a,b,c;h^2)$} \\ \hline 5 &
1080 & 480 & \multicolumn{4}{c|}{} \\ \cline{1-5} 4 & 10728 & 5808 & 2592 &
984 & \multicolumn{2}{c|}{} \\ \hline 3 & 95760 & 61872 & 32988 & 14736 &
5664 & 1920 \\ \hline 2 & 747952 & 575864 & 366256 & 193744 & 86864 &
33832 \\ \hline 1 & 5169728 & 4692096 & 3544608 & 2224088 & 1170192 &
526608 \\ \hline 0 & 33071072 & 34336864 & 30314016 & 22428704 & 13882384
& 7264872 \\ \hline $a\to$ & 0 & 1 & 2 & 3 & 4 & 5 \\ \hline \end{tabular}
\par
\smallskip
\begin{tabular}{|r|c|c|c|c|c|c|} \hline
$c\downarrow$ & \multicolumn{6}{|c|}{$C_5(a,b,c;h^2)$} \\
\hline 2 & 11704 & 3656 & \multicolumn{4}{c|}{} \\ \cline{1-5} 1 & 207356 &
72888 & 23208 & 6752 & \multicolumn{2}{c|}{} \\ \hline 0 & 3276944 &
1300816 & 462744 & 149504 & 44272 & 12024 \\ \hline $a\to$ & 6 & 7 & 8 & 9
& 10 & 11 \\ \hline \end{tabular}
\par
\smallskip
{\bf Table 7.7} Characteristic numbers for plane rational 1-cuspidal quintics;
cusp at a specified point.
\end{center}
\par
\bigskip
Equation $1134$:
For any $a$, $b$ and $c$ with $a+b+2c=3d-3$ we have that
\begin{align*}
C_d(a,b&,c;h)=-bC_d(a,b-1,c;h^2)
+\frac{1}{d^2}\bigg\{d(2d-2)N_d(a,b,c+1)+(2-d)N_d(a+1,b+1,c)\\
+dN_d(a,b+2,c)&+\sum_{\substack{d_1+d_2=d\\
a_1+a_2=a\\b_1+b_2=b+1\\c_1+c_2=c}}
N_{d_1}(a_1,b_1,c_1)E_{d_2}(a_2,b_2,c_2;\hd^2)
\bigg[2d_1d_2^2\binom{b}{b_1-1}
-2d_1^2d_2\binom{b}{b_1}\bigg]\binom{a}{a_1}\binom c{c_1}\\
&+\sum_{\substack{d_1+d_2=d\\
a_1+a_2=a+1\\b_1+b_2=b+1\\c_1+c_2=c}}
N_{d_1}(a_1,b_1,c_1)E_{d_2}(a_2,b_2,c_2;\hd) \bigg[2d_2\binom{b}{b_1-1}
-2d_1^2d_2\binom{b}{b_1}
\bigg]\binom{a}{a_1-1}\binom c{c_1}
\bigg\}
\end{align*}
\par
In the tables of the numbers $C_d(a,b,c;h)$ the number of tangent conditions
$b$ is given by $3d-3-a-2c$.
\par
\begin{center}
\par
\begin{tabular}{|r|c|c|c|c|c|c|c|} \hline $c\downarrow$ &
\multicolumn{7}{|c|}{$C_3(a,b,c;h)$} \\ \hline 3 & 6 &
\multicolumn{6}{c|}{} \\ \cline{1-4} 2 & 20 & 18 & 10 &
\multicolumn{4}{c|}{} \\ \cline{1-6} 1 & 42 & 60 & 54 & 30 & 12 &
\multicolumn{2}{c|}{} \\ \hline 0 & 72 & 132 & 186 & 168 & 96 & 42 & 12 \\
\hline $a\to$ & 0 & 1 & 2 & 3 & 4 & 5 & 6\\ \hline \end{tabular}
\par
\smallskip
{\bf Table 7.8} Characteristic numbers for plane rational 1-cuspidal cubics;
cusp on a specified line.
\end{center}
\par
\bigskip
The numbers for cubics were calculated by Schubert \cite{Schubert}, and
Aluffi \cite{Aluffi2}.
\par
\begin{center}
\par
\begin{tabular}{|r|c|c|c|c|c|c|c|c|c|c|} \hline $c\downarrow$ &
\multicolumn{10}{|c|}{$C_4(a,b,c;h)$} \\ \hline 4 & 210 & 120 &
\multicolumn{8}{c|}{} \\ \cline{1-5} 3 & 1200 & 900 & 510 & 228 &
\multicolumn{6}{c|}{} \\ \cline{1-7} 2 & 5640 & 5404 & 3968 & 2232 & 1006 &
380 & \multicolumn{4}{c|}{} \\ \cline{1-9}
1 & 21844 & 26344 & 24812 & 17956 & 10084 & 4604 & 1774 & 592 &
\multicolumn{2}{c|}{} \\ \hline
0 & 81324 & 111776 & 128992 & 118296 & 84284 & 47284 & 21816 & 8552 & 2926
& 864 \\ \hline
$a\to$ & 0 & 1 & 2 & 3 & 4 & 5 & 6 & 7 & 8 & 9 \\ \hline \end{tabular}
\par
\smallskip
{\bf Table 7.9} Characteristic numbers for plane rational 1-cuspidal quartics;
cusp on a specified line.
\end{center}
\par
\bigskip
\begin{center}
\par
\begin{tabular}{|r|c|c|c|c|c|c|} \hline
$c\downarrow$ & \multicolumn{6}{|c|}{$C_5(a,b,c;h)$} \\ \hline 6 & 1800
& \multicolumn{5}{c|}{} \\ \cline{1-4} 5 & 17160 & 9600 & 17160 &
\multicolumn{3}{c|}{} \\ \cline{1-6} 4 & 144024 & 96864 & 53520 & 24744 &
9840 & \multicolumn{1}{c|}{} \\ \hline 3 & 1068168 & 851832 & 561660 &
307056 & 141864 & 56808 \\ \hline 2 & 6944128 & 6542664 & 5123864 &
3327352 & 1805040 & 834488 \\ \hline 1 & 41321296 & 44794720 & 41150000 &
31591528 & 20223208 & 10897256 \\ \hline 0 & 239546016 & 286715392 &
299655536 & 267383808 & 200794048 & 126634616 \\ \hline
$a\to$ & 0 & 1 & 2 & 3 & 4 & 5 \\ \hline \end{tabular}
\par
\smallskip
\begin{tabular}{|r|c|c|c|c|c|c|c|} \hline $c\downarrow$ &
\multicolumn{7}{|c|}{$C_5(a,b,c;h)$} \\ \hline 3 & 20172 &
\multicolumn{6}{c|}{} \\ \cline{1-4} 2 & 336544 & 120864 & 39272 &
\multicolumn{4}{c|}{} \\ \cline{1-6} 1 & 5040940 & 2045080 & 741368 &
243568 & 73216 &
\multicolumn{2}{c|}{} \\ \hline
0 & 67751352 & 31325416 & 12761768 & 4659408 & 1544416 & 469112 & 130896
\\ \hline
$a\to$ & 6 & 7 & 8 & 9 & 10 & 11 & 12 \\ \hline \end{tabular}
\par
\smallskip
{\bf Table 7.10} Characteristic numbers for plane rational 1-cuspidal quintics;
cusp on a specified line.
\end{center}
\par
\bigskip
Equation $1133$:
For any $a$, $b$ and $c$ with $a+b+2c=3d-2$ we have that
\begin{align*}
C_d(a,b&,c;1)=-bC_d(a,b-1,c;h)-cC_d(a,b,c-1;h^2)-\binom{b}{2}C_d(a,b-
2,c;h^2)\\ &+\frac{1}{d^2}\bigg\{4(d-1)N_d(a+1,b,c)+(3d^2-4d)N_d(a,b+1,c)
\\ &-4\cdot\sum_{\substack{d_1+d_2=d\\
a_1+a_2=a\\b_1+b_2=b\\c_1+c_2=c}}
N_{d_1}(a_1,b_1,c_1)E_{d_2}(a_2,b_2,c_2;\hd^2)
d_1d_2\binom{a}{a_1}\binom{b}{b_1}\binom{c}{c_1}\\ &-
4\cdot\sum_{\substack{d_1+d_2=d\\
a_1+a_2=a+1\\b_1+b_2=b\\c_1+c_2=c}}
N_{d_1}(a_1,b_1,c_1)E_{d_2}(a_2,b_2,c_2;\hd) d_2^2\binom{a}{a_1-
1}\binom{b}{b_1}\binom{c}{c_1} \bigg\}
\end{align*}
\par
In the tables of the numbers $C_d(a,b,c;1)$ the number of tangent conditions
$b$ is given by $3d-2-a-2c$. \begin{center}
\par
\begin{tabular}{|r|c|c|c|c|c|c|c|c|} \hline $c\downarrow$ &
\multicolumn{8}{|c|}{$C_3(a,b,c;1)$} \\ \hline 3 & 6 & 6 &
\multicolumn{6}{c|}{} \\ \cline{1-5} 2 & 12 & 18 & 18 & 12 &
\multicolumn{4}{c|}{} \\ \cline{1-7} 1 & 18 & 36 & 54 & 54 & 36 & 18 &
\multicolumn{2}{c|}{} \\ \hline 0 & 24 & 60 & 114 & 168 & 168 & 114 & 60 &
24\\ \hline $a\to$ & 0 & 1 & 2 & 3 & 4 & 5 & 6 & 7\\ \hline \end{tabular}
\par
\smallskip
{\bf Table 7.11} Characteristic numbers for plane rational 1-cuspidal cubics.
\end{center}
\par
\bigskip
The numbers for cubics were calculated by Schubert \cite{Schubert}, by Kleiman
and Speiser \cite{KleimanSpeiser2}, and by Aluffi \cite{Aluffi2}. Note how
the numbers reflect the fact that the dual of a cuspidal cubic is also a
cuspidal
cubic, with duality between point and tangent conditions. \par
Given a plane cuspidal cubic $C$, under the morphism $C\to
\check{C}$ the image of the cusp of $C$ is the point of flex tangency on
$\check{C}$ and vice versa. Therefore we have the following alternative
enumerative significance: $C_3(a,b,c;h^2)$ (listed in Table 7.5) is equal to
the number of cuspidal
cubics with given flex, through $b$ points, tangent to $a$ lines and
incident to
$c$ flags, and $C_3(a,b,c;h)$ (listed in Table 7.8) is equal to the number of
cuspidal cubics with
flex tangent through a given point, through $b$ points, tangent to $a$ lines
and incident to $c$ flags.
\par
\begin{center}
\par
\begin{tabular}{|r|c|c|c|c|c|c|c|c|c|c|c|} \hline $c\downarrow$ &
\multicolumn{11}{|c|}{$C_4(a,b,c;1)$} \\ \hline 5 & 90 &
\multicolumn{10}{c|}{} \\ \cline{1-4} 4 & 450 & 360 & 216 &
\multicolumn{8}{c|}{} \\ \cline{1-6} 3 & 2016 & 2016 & 1566 & 936 & 450 &
\multicolumn{6}{c|}{} \\ \cline{1-8} 2 & 7344 & 9228 & 9096 & 6948 & 4134 &
2004 & 828 & \multicolumn{4}{c|}{} \\ \cline{1-10}
1 & 24012 & 35568 & 43452 & 42012 & 31644 & 18792 & 9198 & 3852 & 1422 &
\multicolumn{2}{c|}{} \\ \hline
0 & 75924 & 126720 & 180288 & 212976 & 201132 & 149364 & 88560 & 43668 &
18486 & 6912 & 2304\\ \hline
$a\to$ & 0 & 1 & 2 & 3 & 4 & 5 & 6 & 7 & 8 & 9 & 10 \\ \hline \end{tabular}
\par
\smallskip
{\bf Table 7.12} Characteristic numbers for plane rational 1-cuspidal
quartics.
\end{center}
\par
\bigskip
\begin{center}
\par
\begin{tabular}{|r|c|c|c|c|c|c|c|} \hline $c\downarrow$ &
\multicolumn{7}{|c|}{$C_5(a,b,c;1)$} \\ \hline 6 & 7560 & 4320 &
\multicolumn{5}{c|}{} \\ \cline{1-5} 5 & 58680 & 41040 & 23400 & 11160 &
\multicolumn{3}{c|}{} \\ \cline{1-7} 4 & 417528 & 343584 & 234432 & 132408
& 63216 & 26208 & \multicolumn{1}{c|}{} \\ \hline
3 & 2600136 & 2522808 & 2038572 & 1367856 & 766296 & 365976 & 152748 \\
\hline
2 & 14519232 & 16336416 & 15523656 & 12319632 & 8150736 & 4535952 &
2167776 \\ \hline
1 & 76377456 & 96562656 & 105583536 & 97988544 & 76305384 & 49814136 &
27552636 \\ \hline
0 & 392798880 & 543805632 & 664607952 & 704860128 & 637644672 & 486730080
& 313485192 \\ \hline
$a\to$ & 0 & 1 & 2 & 3 & 4 & 5 & 6 \\ \hline \end{tabular}
\end{center}
\par
\smallskip
\begin{center}
\begin{tabular}{|r|c|c|c|c|c|c|c|} \hline $c\downarrow$ &
\multicolumn{7}{|c|}{$C_5(a,b,c;1)$} \\ \hline 3 & 56952 &
\multicolumn{6}{c|}{} \\ \cline{1-4} 2 & 910056 & 342312 & 117216 &
\multicolumn{4}{c|}{} \\ \cline{1-6} 1 & 13170816 & 5553864 & 2103264 &
725616 & 230616 &
\multicolumn{2}{c|}{} \\ \hline
0 & 172272672 & 82282248 & 34794432 & 13240080 & 4592952 & 1467792 &
435168
\\ \hline
$a\to$ & 7 & 8 & 9 & 10 & 11 & 12 & 13\\ \hline \end{tabular}
\par
\smallskip
{\bf Table 7.13} Characteristic numbers for plane rational 1-cuspidal
quintics.
\end{center}
\par
\bigskip
The four series of characteristic numbers of Theorem \ref{algthm} have been
calculated up to degree 10 and are available upon request to the authors. The
Maple source code for the calculation may also be requested.
\newpage
\ifx\undefined\bysame
\newcommand{\bysame}{\leavevmode\hbox to3em{\hrulefill}\,}
\fi


\begin{thebibliography}{Kon95}

\bibitem[Alu91]{Aluffi2}
P.~Aluffi, {\em The enumerative geometry of plane cubics. II. Nodal and
  cuspidal cubics}, Math. Annalen. {\bf 289} (1991), 543--572.

\bibitem[Beh96]{Behrend}
K.~Behrend, {\em Gromov-Wittwn invariants in algebraic geometry}, 1996,
  alg-geom/9601011.

\bibitem[BM95]{BehrendManin}
K.~Behrend and Y.~Manin, {\em Stacks of stable maps and Gromov-Witten
  invariants}, 1995, alg-geom/9506023.

\bibitem[CK91]{CKtrip}
S.~J. Colley and G.~Kennedy, {\em Triple and quadruple contact of plane
  curves}, Enumerative Algebraic Geometry (Proc. 1989 Zeuthen Symposium) (S.~L.
  Kleiman and A.~Thorup, eds.), Contemporary Math., vol. 123, Amer. Math. Soc.,
  Providence, R.I., 1991, pp.~31--59.

\bibitem[Col88]{Collino}
A.~Collino, {\em Evidence for a conjecture of Ellingsrud and Str\o mme on the
  Chow ring of $\operatorname{Hilb}_d(\pp)$}, Illinois J. Math. {\bf 32}
  (1988), 171--210.

\bibitem[EG96]{EdidinGraham}
D.~Edidin and W.~Graham, {\em Equivariant intersection theory}, Preprint, 1996.

\bibitem[FI95]{FrancescoItzykson}
P.~Di Francesco and C.~Itzykson, {\em Quantum Intersection rings}, The Moduli
  Space of Curves, Birkh{\"a}user, Boston, 1995, pp.~81--148.

\bibitem[Ful95]{Fulton2}
W.~Fulton, {\em Enumerative Geometry via Quantum Cohomology}, Lecture notes AMS
  summer institute Santa Cruz, 1995.

\bibitem[Kle74]{Kleiman4}
S.~L. Kleiman, {\em The transversality of a general translate}, Compositio
  Math. {\bf 28} (1974), no.~3, 287--297.

\bibitem[KM94]{KontsevichManin}
M.~Kontsevich and Y.~Manin, {\em Gromov-Witten classes, quantum cohomology and
  enumerative geometry}, 1994, hep-th/9402147.

\bibitem[Kon95]{Kontsevich}
M.~Kontsevich, {\em Enumeration of rational curves via torus action}, The
  Moduli Space of Curves, Birkh{\"a}user, Boston, 1995, pp.~335--368.

\bibitem[KS86]{KleimanSpeiser2}
S.~L. Kleiman and R.~Speiser, {\em Enumerative geometry of cuspidal plane
  cubics}, (Proc. Alg. Geom. Vancouver 1984), Canad. Math. Soc. Conf. Proc.,
  vol.~6, AMS, Providence, R.I., 1986, pp.~227--267.

\bibitem[KS87]{KleimanSpeiser}
S.~L. Kleiman and R.~Speiser, {\em Enumerative geometry of nodal plane curves},
  (Algebraic Geometry, Sundance 1986), Lecture Notes in Math., vol. 1311,
  Springer-Verlag, Berlin, 1987.

\bibitem[LT95]{LiTian}
J.~Li and G.~Tian, {\em The quantum cohomology of homogeneous varieties}, 1995,
  alg-geom/9504009.

\bibitem[LT96]{LiTian2}
J.~Li and G.~Tian, {\em Virtual moduli cycles and GW-invariants}, 1996,
  alg-geom/9602007.

\bibitem[Mai71]{Maillard}
S.~Maillard, {\em Recherche des caract\'eristiques des syst\`emes
  \'el\'ementaires de courbes planes du troisi\`eme ordre}, 1871, Th\`eses
  pre\'esent\'ees \`a la Facult\'e des Sciences de Paris.

\bibitem[Pan94]{Pandharipande}
R.~Pandharipande, {\em `Notes on Kontsevich's compactification of the space of
  maps}, preprint, 1994.

\bibitem[Pan95]{Pandharipande2}
R.~Pandharipande, {\em Intersections of $\q$-divisors on Kontsevich's moduli
  space ${\overline M}^0_{0,n}(\bp^r,d)$ and enumerative geometry}, 1995,
  alg-geom/9504004.

\bibitem[Sac85]{Sacchiero}
G.~Sacchiero, {\em Numeri caratteristici delle cubiche piane cuspidali}, 1985,
  preprint.

\bibitem[Sch79]{Schubert}
H.~C.~H. Schubert, {\em Kalk{\"u}l der abz{\"a}hlende geometrie}, B. G.
  Teubner, Leipzig, 1879, reprinted by Springer-Verlag 1979 with an
  introduction by S. L. Kleiman.

\bibitem[Vis89]{Vistoli}
A.~Vistoli, {\em Intersecion theory on algebraic stacks and their moduli}, Inv.
  Math. {\bf 97} (1989), 613--670.

\bibitem[Zeu72]{Zeuthen}
H.~Zeuthen, {\em D\'etermination des charact\'eristiques des syst\`emes
  \'el\'ementaires de cubiques}, C. R. Acad. Sci. Paris {\bf 74} (1872),
  521--526, 604--607, 726--729.

\end{thebibliography}
\end{document}